\documentclass[11pt,a4paper]{article}
\usepackage{authblk}
\usepackage{graphicx}
\usepackage[bottom]{footmisc}
\usepackage{ulem} 
\usepackage{newtxtext} 
\usepackage[varvw]{newtxmath} 
\usepackage{tikz}
\usepackage[hidelinks]{hyperref}
\usepackage{amsmath}
\usepackage{array}
\usepackage{bbm}
\usepackage{empheq}
\usepackage{enumitem}
\usepackage{cite}
\usepackage[left=10mm, right=10mm, top=20mm, bottom=20mm, heightrounded, marginparsep=5mm]{geometry} 
\usepackage{caption}
\usetikzlibrary{decorations.markings}

\newcommand{\spt}{\mathop\mathrm{spt}}
\newcommand{\R}{\mathbb{R}}
\newcommand{\N}{\mathbb{N}}
\newcommand{\tv}{\mathrm{TV}}
\newcommand{\C}[1]{\mathbf{C^{\pmb #1}}}
\newcommand{\Lip}{\mathbf{Lip}}
\newcommand{\Cc}[1]{\mathbf{C_c^{#1}}}
\renewcommand{\L}[1]{\mathbf{L^{\pmb #1}}}
\newcommand{\W}[1]{\mathbf{W^{\pmb #1}}}
\newcommand{\BV}{\mathbf{BV}}
\newcommand{\Lloc}[1]{\mathbf{L^{#1}_{loc}}}
\renewcommand{\d}{\mathrm{d}}
\renewcommand{\epsilon}{\varepsilon}

\renewcommand{\geq}{\geqslant}
\renewcommand{\leq}{\leqslant}

\usepackage{xcolor}     
\usepackage{lipsum}     
\usepackage{ntheorem}   
\usepackage{mdframed}   
\theoremstyle{break}
\theoremheaderfont{\bfseries}
\newmdtheoremenv[%
linecolor=gray,
backgroundcolor=gray!40,%
innertopmargin=5pt,%
innerbottommargin=5pt,%
innerleftmargin=5pt,%
innerrightmargin=5pt,%
skipbelow=0pt,%
skipabove=4pt,%
ntheorem]{svgraybox}{Hypothesis}[section]

\newtheorem{proposition}{Proposition}[section]
\newtheorem{theorem}[proposition]{Theorem}

\newtheorem{definition}[proposition]{Definition}

\newtheorem{remark}[proposition]{Remark}

\begin{document}
\normalem
\allowdisplaybreaks

\title{The mathematical theory of Hughes' model: a survey of results}
\author[1]{D.~Amadori}
\affil[1]{\footnotesize Department of Information Engineering, Computer Science and Mathematics, University of L'Aquila 67100, Italy} 
\author[2,3]{B.~Andreianov}
\affil[2]{Institut Denis Poisson (CNRS UMR7013), Universit{\'e} de Tours, Universit{\'e} d'Orl{\'e}ans, Parc Grandmont 37200 Tours, France}
\affil[3]{Peoples' Friendship University of Russia (RUDN University) 6 Miklukho-Maklaya St, Moscow, 117198, Russian Federation}
\author[1]{M.~Di~Francesco}
\author[1]{S.~Fagioli}
\author[2]{T.~Girard}
\author[4]{P.~Goatin}
\affil[4]{Universit\'e C\^ote d'Azur, Inria, CNRS, LJAD, Sophia Antipolis, France}
\author[5,6]{P.~Markowich}
\affil[5]{Mathematical and Computer Sciences and Engineering Division, King Abdullah University of Science and Technology, Thuwal 23955-6900, Kingdom of Saudi Arabia}
\affil[6]{Faculty of Mathematics, University of Vienna, Oskar-Morgenstern-Platz 1, 1090 Vienna, Austria} 
\author[7]{J.-F.~Pietschmann}
\affil[7]{Institute of Mathematics, University of Augsburg, Universitätsstr. 12a, 86159 Augsburg, Germany} 
\author[8,9,10]{M.D.~Rosini}
\affil[8]{Department of Mathematics and Computer Science, University of Ferrara 44121, Italy}
\affil[9]{Department of Management and Business Administration, University “G.~d'Annunzio” of Chieti-Pescara 65127, Italy}
\affil[10]{Uniwersytet Marii Curie-Sklodowskiej, Plac Marii Curie-Sklodowskiej 1 20-031 Lublin, Poland} 
\author[11]{G.~Russo}
\affil[11]{Department of Mathematics and Computer Science, University of Catania, Viale Andrea Doria 6, 95125, Catania, Italy}
\author[1]{G.~Stivaletta}
\author[12]{M.T.~Wolfram}
\affil[12]{Mathematics Institute, University of Warwick, Coventry CV4 7AL, United Kingdom} 

\maketitle
\normalsize
\abstract{
We provide an overview of the results on Hughes' model for pedestrian movements  available in the literature. The model consists of a nonlinear conservation law coupled with an eikonal equation. The main difficulty in developing a proper mathematical theory lies in the 
lack of regularity of the flux in the conservation law, which yields the possibility of non-classical shocks that are generated non-locally by the whole distribution of pedestrians. This is a possible reason behind the availability of existence results only on one-dimensional spatial domains, despite the model having a more natural setting in two spatial dimensions.
\\
After the first successful approaches to solving a regularised version of the model, researchers focused on the structure of the Riemann problem, which led to local-in-time existence results for Riemann-type data and paved the way for a WFT (\textit{Wave-Front Tracking}) approach to the solution semigroup. In parallel, a DPA (\textit{Deterministic Particles Approximation}) approach was developed in the spirit of follow-the-leader approximation results for scalar conservation laws. Beyond having proved to be powerful analytical tools, the WFT and
the DPA approaches also led to interesting numerical results.
\\
However, only existence theorems on very specific classes of initial data (essentially ruling out non-classical shocks) have been available until very recently. A proper existence result using a DPA approach was proven not long ago in the case of a linear coupling with the density in the eikonal equation. Shortly after, a similar result was proven via a fixed point approach.
\\
We provide a detailed statement of the aforementioned  results and sketch the main proofs. We also provide a brief overview of results that are related to Hughes' model, such as the derivation of a dynamic version of the model via a mean-field game strategy, an alternative optimal control approach, and a localized version of the model. We also present the main numerical results within the WFT and DPA frameworks.
}

\section{Introduction}

In recent years, the flow of pedestrians has attracted remarkable scientific interest due to its potential in multidisciplinary applications, for instance in the design of safety systems in case of evacuation of a building (such as a stadium) or in crowd management during events of mass gathering.

Three main different types of modelling approaches are present in the literature, corresponding to three different levels of resolution:
\begin{itemize}
    \item a description of the state of the individuals and of the interactions among them, referred to as the \emph{individual}, or \emph{microscopic}, level;
    \item a statistical description of a sample of the system, known as the \emph{mesoscopic}, or \emph{kinetic}, level;
    \item a continuum description at the level of the interaction of sub-populations of the system, known as the \emph{macroscopic} level.
\end{itemize}
Each level has an associated class of mathematical equations which provide an appropriate model. 
Usually their structures at various levels are completely different, see e.g.~\cite{MR4059365}.
We defer the reader to \cite{MR4180806, MR3932134, bellomocrowd} for an overview of the various research directions in the field of crowd dynamics.

One of the most original and mathematically challenging models is the one proposed by Roger L. Hughes~\cite{Hughes02}, which is part of the third approach above.
Hughes' model describes evacuation scenarios, in which a crowd wants to exit a given domain $\mathsf{D} \subset \R^2$, with one or several exits, as fast as possible. The crowd population is described through a density $\rho$. 
The driving force towards the exits is the gradient of a potential $\phi$, which satisfies an eikonal equation coupled with the density $\rho$.
The potential $\phi$ represents the expected travel time towards an exit and accounts for the best strategy to minimise the exit time. The resulting model is a nonlinear conservation law for $\rho$ coupled with the gradient of the potential $\phi$; the latter depends on $\rho$ non-locally in space.

Even in the simplest case of space dimension one, no more than Lipschitz continuity can be expected for $\phi$.
Moreover, $\phi_x$ can change its sign just once from positive to negative.
In this case, a turning curve $x=\xi(t)$ may be defined in the one-dimensional domain, at which $\phi$ reaches its maximum. As a result, the flux of the conservation law for $\rho$ is possibly discontinuous along $x=\xi(t)$ and, on the other hand, $\xi$ depends non-locally on $\rho$.
Furthermore, $\rho$ is not expected to satisfy in general the Lax entropy inequalities (see \cite{MR1927887}). Therefore, the possible appearance of non-classical shocks along the turning curve has to be taken into account.
Non-classical shocks correspond to pedestrians changing direction during the evacuation.
Because of these various difficulties, the existence and uniqueness analysis for the Hughes’ model appears to be challenging.
This motivates the first approach to the problem developed in \cite{markowich}, in which a \emph{smoothened version} of the eikonal equation is considered, with an extra Laplacian term. Despite it covers only the one-dimensional case and despite it deals only with an approximated version of the model rather than the actual model, the result in \cite{markowich} remains until now the only result in a large data setting which holds for a fairly general class of density-potential coupling. As we will detail later on, the first existence results for the actual Hughes model appeared only very recently, and only for specific couplings.
Let us also mention the existence and uniqueness result proved in \cite[Theorem~2.6]{MR3823842} for a two dimensional regularized version of the Hughes model, which holds for large data but with a potential $\phi$ depending only on the given domain $\mathsf{D}$, see also \cite{MR4096596}.

Parallel to \cite{markowich} or shortly after it, some researchers started to study the Riemann problem for the model in the spirit of scalar conservation laws and to develop proper numerical schemes, see \cite{AmadoriDiFrancesco, AmadoriGoatinRosini, ElKhatibGoatinRosini, AndRosSti} for the Riemann data part and \cite{MR3149318, MR3177723, MR3619091, MR3698447, MR3055243, HUANG2009127, MR3277564} for the numerics. Both in \cite{AmadoriDiFrancesco} and \cite{ElKhatibGoatinRosini}, the authors study the Riemann problems for the Hughes' model in detail.
This study is strictly related to the effectiveness of the Wave-Front Tracking (WFT) strategy \cite{DafermosWFT} for the Hughes' model.
The WFT algorithm was then first exploited in \cite{MR3055243}, but only for numerical purposes, and then in \cite{AmadoriGoatinRosini} to prove the first existence result for the Hughes' model, but under very restrictive assumptions that rule out non-classical shocks.
A simpler proof of an analogous existence result was then obtained in \cite{DiFrancescoFagioliRosiniRussoKRM} by means of a Deterministic Particle Approximation (DPA) and the results proved in \cite{DiFrancescoRosini}, see also \cite{DiFrancescoFagioliRosini-BUMI, MR3644595, DiFrancescoFagioliRosiniRusso, MR4026959}.

The first existence result accounting for the possible presence of non-classical shocks was recently obtained in \cite{AndRosSti}.
The authors obtain this result by exploiting the properties of the linear cost introduced in \cite{ElKhatibGoatinRosini} (the key fact here is that linear costs yield a uniform Lipschitz bound on $\xi$), combined with the DPA adapted to the Hughes' model in \cite{DiFrancescoFagioliRosiniRussoKRM}. 
Despite being only valid for linear costs and in one space dimension, this result has the merit of being the \emph{first} existence result on the Hughes model for large data and in presence of non-classical shocks. This result is re-proved in \cite{AndrGirard-preprint} via a non-constructive Schauder fixed-point approach allowing for a wide variety of generalizations of the one-dimensional Hughes' model (different ways to compute the turning curve $\xi$ from the density $\rho$, different exit conditions). 

Concerning uniqueness, only very partial results are available for the one-dimensional Hughes' model; they require  $\BV$ regularity of the density $\rho$ and the highly restrictive assumption of zero density traces $\rho(t,\xi(t)^\pm)$ at the turning curve (see \cite[Theorem~4]{AndRosSti}, see also \cite{AmadoriGoatinRosini,DiFrancescoFagioliRosiniRussoKRM} for particular cases).

Apart from the regularised version proposed in \cite{markowich}, other variants of the Hughes model have been proposed: a first one in \cite{burger, HerzogPietschmannWinkler} obtained a similar model with a time derivative in the eikonal equation, justified through an optimal control problem, and a second one in \cite{carrillo} trying to remove the \emph{global awareness} of the pedestrians in the model, which seems unrealistic in some situations. Further variants with more flexible boundary conditions for the density, with memory or relaxation effects in the dynamics of $\xi$, are proposed and studied in \cite{AndrGirard-preprint}.

The chapter is structured as follows. 
In Section~\ref{sec:construction} we derive Hughes' model in the way it was done in the original paper~\cite{Hughes02} by Roger L.~Hughes, plus some additional considerations by the authors of this survey. 
We also provide a rephrasement of the model in the special case of one space dimension. 
In Section~\ref{sec:RP} we detail the local-in-time solution of the Riemann problem. 
In Section~\ref{sec:existence} we collect the existence result provided for the model, from the ones holding only for small data or symmetric data provided in \cite{AmadoriGoatinRosini,DiFrancescoFagioliRosiniRusso}, to the main one provided in \cite{AndRosSti} for the case of linear cost. 
In Section~\ref{s:WFT} we describe the construction of the Wave-Front Tracking (WFT) algorithm used to prove the existence results in \cite{AmadoriGoatinRosini}. 
In Section~\ref{s:mpa} we introduce the Deterministic Particle Approximation (DPA) of the model leading to the results in \cite{DiFrancescoFagioliRosiniRusso} and \cite{AndRosSti}. 
In Section~\ref{s:linear} we describe in detail the main existence result of \cite{AndRosSti}.
In Section~\ref{s:fixed-point} we briefly describe the fixed-point approach of \cite{AndrGirard-preprint}, with a second proof of this main existence result and several extensions.
Section~\ref{sec:sim} is devoted to numerical simulations, both using the WFT algorithm and the DPA scheme. 
Finally, in Section~\ref{sec:modified} we summarise the modified versions of the model considered in this survey, namely the smoothed version proposed in \cite{markowich}, the dynamic one of \cite{burger}, optimal control problems in \cite{HerzogPietschmannWinkler}, and the localised one of \cite{carrillo}.

This survey paper covers a very high percentage of the available work on Hughes' model. The model is, however, becoming quite well-known in the applied mathematics community, and it is therefore quite likely that we may have missed some papers. The main focus of the present paper is the well-posedness theory, for which we believe we covered the main results available in the literature. Our choice of the numerical results is intentional: we cover the WFT and the DPA approaches because they are strictly related to the techniques used in some of the existence results described here. Our choice of the extended models is also not accidental: apart from the result on the regularized model, we addressed results which are relevant in that they feature slight modifications of the model which are significant from the point of view of the applications. 

\section{Construction of the model}\label{sec:construction}

\subsection{The two-dimensional case}

Back in 2002, Roger L.~Hughes proposed a model for a two-dimensional flow of pedestrians~\cite{Hughes02}.
The model accounts for the possible presence of obstacles (walls, columns, etc.) and for multiple pedestrian types. 
For simplicity, in this review we consider the case in which only one population of pedestrians is involved.
Since the movement of pedestrians takes place in a two dimensional space, the model is typically set on a bounded domain $\mathsf{D}\in \R^2$.
The pedestrian flow is described in terms of two quantities:
\begin{itemize}
\item
\emph{density}, $\rho=\rho(t,X)$, which is the number of individuals per unit area at a given time $t$ and location $X=(x_1,x_2)\in \mathsf{D}$ of the walking space, and
\item
\emph{velocity}, $V = V(t,X) = \left(V_1(t,X),V_2(t,X)\right)\in \R^2$, which is the average velocity of individuals located within a unit area of the walking space, at a given time $t$ and location $X=(x_1,x_2)\in \mathsf{D}$.
\end{itemize}
The conservation of the number of pedestrians is expressed by the continuity equation
\begin{equation}\label{e:aHughes2D}
\rho_t + (\rho V_1)_{x_1} + (\rho V_2)_{x_2} = 0.
\end{equation}
Recall that the above equation is obtained by equating the net flow of pedestrians into a small region to the rate of accumulation of pedestrians in the region, and then letting the area of the region shrink to zero, see for instance \cite{Dafermosbook}.

To complete the model, we assume what follows:

\begin{svgraybox}
The speed of pedestrians is determined as a function of the density $v=v(\rho)$, with $v:[0,\rho_{\max}] \rightarrow [0,v_{\max}]$ being a decreasing function such that $v(0) = v_{\max}>0$ and $v(\rho_{\max}) = 0$.
\end{svgraybox}

\begin{svgraybox}
Each pedestrian has a common sense of the task (called \emph{potential}) they face to reach their common destination.
In particular, two individuals at different locations having the same potential don't see any advantage in exchanging their positions.
Furthermore, pedestrians move orthogonally to level set curves of the potential.
\end{svgraybox}

\begin{svgraybox}
Pedestrians seek the path that minimizes their given travel cost.
\end{svgraybox}

The first hypothesis is standard in traffic modeling, as lower speeds correspond to higher densities, see for instance \cite{GaravelloPiccolibook, Rosinibook}.
The parameter $\rho_{\max}$ stands for the maximum admissible density and $v_{\max}$ for the maximum speed.
The second hypothesis implies that each pedestrian knows the overall density distribution of the crowd, that occurs if, for instance, shorter pedestrians take their direction from the tallest pedestrians, who have an overall view of the situation.
The third hypothesis is about travel cost as discussed later in this section.

By the first hypothesis, the velocity components are given by
\begin{align}
\label{e:velocity2D}
V_1&=v(\rho) \hat{\phi}_{1},&
V_2&=v(\rho) \hat{\phi}_{2},
\end{align}
where $\hat{\phi}_{1}$ and $\hat{\phi}_{2}$ are the direction cosines of the motion and $v(\rho)$ is the speed.

A consequence of the second hypothesis is that there is no advantage of moving along a line of constant potential, but that actually pedestrians move down the gradient of the potential $\phi$, that is, their trajectories are parallel to the gradient of $\phi$:
\begin{equation}
\label{e:directions2D}
\left(\begin{array}{@{}c@{}}
\hat{\phi}_{1}
\\
\hat{\phi}_{2}
\end{array}\right)
=
-\frac{1}{\|\nabla\phi\|} \nabla\phi.
\end{equation}

Assume $X_1, X_2\in \mathsf{D}$ are very close and with their connecting line perpendicular to the level sets $\phi(X)=\phi(X_1)$ and $\phi(X)=\phi(X_2)$. 
If $\phi(X_1) > \phi(X_2)$, then
\[\phi(X_1)-\phi(X_2)\simeq\|\nabla\phi\| \ \|X_1-X_2\|.\]
Assuming a pedestrian moves from $X_1$ to $X_2$ with constant speed $\|V\|$ during the time interval $[t_1,t_2]$ gives
\begin{equation}\label{e:marco1}
 \phi(X_1)-\phi(X_2)\simeq (t_2-t_1)\ \|\nabla\phi\| \ \|V\| = (t_2-t_1)\ \|\nabla\phi\|\ v(\rho) .
\end{equation}

Recall that, due to third hypothesis, the potential itself measures the \lq estimated\rq\ time, in a tempered way that takes into account the density. 
In order to encode such hypothesis in the model, we locally assume that \lq small\rq\ differences in the potential are proportional to the product of the speed and the density, that is
\begin{equation}\label{e:marco2}
 \phi(X_1)-\phi(X_2)\simeq \frac{t_2-t_1}{g(\rho)}, 
\end{equation}
where the factor $g(\rho)$ allows for discomfort at very high densities. 
The function $g \colon [0,+\infty)\rightarrow[1,+\infty)$ is assumed to satisfy $g(0)=1$ (that is, no discomfort when there are no pedestrians around) and to be increasing in $\rho$. Equating \eqref{e:marco1} and \eqref{e:marco2} gives
\begin{equation}
\label{e:aeikonal2D}
\frac{1}{\|\nabla\phi\|} = g(\rho) v(\rho).
\end{equation}

The governing equations are obtained by combining \eqref{e:aHughes2D}, \eqref{e:velocity2D}, \eqref{e:directions2D}, \eqref{e:aeikonal2D} and write
\begin{equation}
\label{e:Hughes2D}
\left\{\begin{array}{@{}>{\displaystyle}l@{}}
\rho_t - \nabla \cdot \left(\rho v(\rho) \frac{\nabla\phi}{\|\nabla\phi\|} \right) = 0,\\
\|\nabla\phi\| = \frac{1}{g(\rho) v(\rho)}.
\end{array}\right.
\end{equation}
The resulting model is therefore a scalar conservation law coupled with an eikonal equation. 
This is the original model formulated by Hughes in \cite{Hughes02}. Later versions of this model \cite{ElKhatibGoatinRosini} regarded the right-hand side of the eikonal equation in \eqref{e:Hughes2D} as a \emph{running cost} $c(\rho)$, which satisfies $c(\rho)=1/\bigl(g(\rho) v(\rho)\bigr)$ in the original formulation by Hughes. 
Typical assumptions on $c$ and $v$ are the following:
\begin{enumerate}
\item[{\bf (H1)}]
The cost map $c \colon [0,\rho_{\max}] \to [1,+\infty)$ is $\C2$, increasing, with $c(0)=1$ and $c''(\rho)\geq0$ for all $\rho\in[0,\rho_{\max}]$. 
\item[{\bf (H2)}]
The speed map $v \colon [0,\rho_{\max}] \to [0,v_{\max}]$ is $\C1$, strictly decreasing, with $v(0)=v_{\max}>0$ and $v(\rho_{\max})=0$.
Moreover, there exists a $\hat{\rho} \in (0,\rho_{\max})$ such that $\left(v(\rho)+\rho v'(\rho)\right) \left(\hat{\rho}-\rho\right) > 0$ for all $\rho \in (0,\rho_{\max}) \setminus \{\hat{\rho}\}$.
\end{enumerate}
Note that $\hat{\rho}$ is the maximum point of $f(\rho) \doteq \rho v(\rho)$.
Below we shall enforce these assumptions.

Model \eqref{e:Hughes2D} requires the specification of an initial condition
\begin{align}
\label{e:HughesIni}
&\rho(0,X)=\bar{\rho}(X),&
&X\in\mathsf{D}.
\end{align}
Typical boundary conditions on $\partial\mathsf{D}$ are characterized by the presence of walls, $\Gamma_w$, corners, $\Gamma_c$, and exits, $\Gamma_e$, so that $\partial\mathsf{D} = \Gamma_w \cup \Gamma_c \cup \Gamma_e$.
We assume pedestrians cannot pass through walls but can move along them.
Hence we impose free-slip boundary conditions, namely $V\cdot\nu=0$ along $\Gamma_w$, where $\nu=\nu(X)$ is the outward unit normal to $\partial\mathsf{D}$ at $X\in \partial\mathsf{D}$.
Hence, model~\eqref{e:Hughes2D} needs the specification of the boundary conditions
\begin{subequations}\label{e:HughesBoundary}
\begin{empheq}[left=\empheqlbrace]{align}
\label{e:HughesBoundaryWp}
&\nabla\phi(t,X) \cdot \nu(X)=0,&
&X\in\Gamma_w,
\\
\label{e:HughesBoundaryE}
&\phi(t,X)=0,&
&X\in\Gamma_e,
\\
\label{e:HughesBoundaryO}
&\rho(t,X)=0,&
&X\in\Gamma_e.
\end{empheq}
\end{subequations}

The boundary condition \eqref{e:HughesBoundaryO} has to be understood in the sense of Bardos, Le Roux and Nedelec \cite{BardosLerouxNedelec}, so it states that for almost every $(t,X) \in (0,+\infty) \times \Gamma_e$ we have
\[\Bigl( \rho(t,X) \, v\bigl(\rho(t,X)\bigr) - k \, v(k) \Bigr) \ \frac{\nabla\phi(t,X)}{\|\nabla\phi(t,X)\|} \cdot \nu(X) \geq 0\]
for all $k \in [0,\rho(t,X)]$.
Intuitively, such a boundary condition is \emph{set valued}.
Indeed, \emph{if} for instance $\nabla\phi(t,X) \cdot \nu(X) \geq 0$, then the above inequality is satisfied if and only if $\rho(t,X)$ belongs to $[0,\hat{\rho}]$, where $\hat{\rho}$ is introduced in {\bf (H2)}.
As we will see, this is the case in the one dimensional setting, see \eqref{NEW:BC_Sect2.2}.
As a consequence, a positive outflow from $\mathsf{D}$ is allowed through the exits thanks to the boundary conditions \eqref{e:HughesBoundaryE} and \eqref{e:HughesBoundaryO}, but not through the walls by \eqref{e:HughesBoundaryWp} and \eqref{e:directions2D}.

\subsection{The one-dimensional case}

As a model for pedestrian movements, Hughes' model natural setting is a two-dimensional space domain. 
However, the mathematical difficulties hidden in the coupling of the conservation law with the eikonal equation in \eqref{e:Hughes2D} led part of the mathematical community (including all the authors of the present manuscript) to focus on its one-dimensional version. 

It turns out that, on a one-dimensional domain, the model can be formulated in an alternative way. 
For simplicity, we pose the model on the one-dimensional interval $\mathsf{C} \doteq (-1,1)$ and assume that two exits are located at $x=\pm1$, so that $\mathsf{C}$ represents, for instance, a corridor or a bridge. 
We therefore obtain
\begin{subequations}
\nonumber\label{e:Hughes1D}
\begin{empheq}[left=\empheqlbrace]{align}
&\rho_t - \left(f(\rho) \frac{\phi_x}{|\phi_x|} \right)_x = 0,
\label{e:Hughes1D-a}
\\
&|\phi_x| = c(\rho),
\label{e:Hughes1D-b}
\end{empheq}
\end{subequations}
with the boundary conditions
\begin{subequations}
\label{e:Hughes1DBoundary}
\begin{empheq}[left=\empheqlbrace]{align}
&\rho(t,-1)=\rho(t,1)=0,&
&t\geq 0,\label{e:Hughes1DBoundary-a}
\\
&\phi(t,-1)=\phi(t,1)=0,&
&t\geq 0.
\label{e:Hughes1DBoundary-b}
\end{empheq}
\end{subequations}
 
Assuming $\rho$ is known in the eikonal equation in \eqref{e:Hughes1D-b}, and assuming as well that $c(\rho)$ is bounded, the eikonal equation can be solved in a \emph{viscosity solution sense}, that is by imposing that $\phi$ is semi-concave, i.e.\ with a second derivative that is bounded from above in the sense of distributions. 
We can therefore assume that $\phi$ is increasing near the boundary point $x=-1$, decreasing near $x=1$, and has a global maximum in $\mathsf{C}$. 
This and the boundary conditions \eqref{e:Hughes1DBoundary-b} for $\phi$ imply
\begin{align*}
 &\phi(t,x)=\int_{-1}^x c\bigl(\rho(t,y)\bigr) \,\d y&&\hbox{for $x$ near $-1$},\\
 &\phi(t,x)=\int_x^1 c\bigl(\rho(t,y)\bigr) \,\d y&&\hbox{for $x$ near $1$}.
\end{align*}
Fix $t>0$.
Since $\phi_x(t,\cdot\,)$ is always non zero, $\phi_x(t,\cdot\,)$ can only have one discontinuity in order to have $\phi(t,\cdot\,)$ in the class of semi-concave functions. 
We call $\xi(t)\in \mathsf{C}$ the discontinuity point. Moreover, to preserve the continuity of $\phi$, we must have
\begin{equation}\label{e:xi1}
 \phi\bigl(t,\xi(t)^-\bigr)=\int_{-1}^{\xi(t)} c\bigl(\rho(t,y)\bigr) \,\d y= \int_{\xi(t)}^1 c\bigl(\rho(t,y)\bigr) \,\d y= \phi\bigl(t,\xi(t)^+\bigr).
\end{equation}
The above calculations give an explicit formula for $\phi$, given the \lq moving discontinuity curve\rq\ $t \mapsto \xi(t)$. 
On the other hand, $\xi(t)$ depends on $\rho(t,\cdot\,)$ via \eqref{e:xi1}. Moreover, for $x\in (-1,\xi(t))$ we have $\phi_x/|\phi_x|=1$ and for $x\in (\xi(t),1)$ we have $\phi_x/|\phi_x|=-1$. 

Therefore, the whole model can be reformulated as follows:
\begin{subequations}
\label{e:Hughes1D_reformulated}
\begin{empheq}[left=\empheqlbrace]{align}
&\rho_t+F(t,x,\rho,\xi)_{x} = 0, &\quad t>0, \ x\in \mathsf{C},
\label{e:Hughes1D_reformulated-a}
\\
&\displaystyle{\int_{-1}^{\xi(t)} c\bigl(\rho(t,y)\bigr) \,\d y= \int_{\xi(t)}^1 c\bigl(\rho(t,y)\bigr) \,\d y,} & \quad t>0, 
\label{e:Hughes1D_reformulated-b}
\\
&\rho(0,x)=\bar\rho(x), & \quad x\in \mathsf{C},
\label{e:Hughes1D_reformulated-c}
\end{empheq}
\end{subequations}
where $\bar\rho$ is the initial datum and in \eqref{e:Hughes1D_reformulated-a} we set
\begin{equation}\label{e:F}
F(t,x,\rho,\xi) \doteq \operatorname{sign}\bigl(x-\xi(t)\bigr) f(\rho),
\end{equation}
coupled with the boundary conditions \eqref{e:Hughes1DBoundary}.
Note that the strong traces of the solution at the boundary points exist due to the genuine nonlinearity of the flux \cite{MR1869441,Panov_traces2} and must satisfy
\begin{align}\label{NEW:BC_Sect2.2}
&\rho(t,-1^+)\leq\hat{\rho},&
&\rho(t,1^-)\leq\hat{\rho},
\end{align}
where $\hat{\rho}$ is introduced in {\bf (H2)} and is the maximum point of $f$.

In principle, the one-dimensional reformulation of the Hughes model \eqref{e:Hughes2D}, \eqref{e:HughesIni}, \eqref{e:HughesBoundary} is represented by \eqref{e:Hughes1D_reformulated} coupled with the Dirichlet boundary conditions \eqref{e:Hughes1DBoundary-a}.
However, in \cite[page~220]{DiFrancescoFagioliRosiniRussoKRM} it is argued that in the one-dimensional case no boundary conditions have to be prescribed.
More rigorously, in \cite[Section~3]{AndRosSti} it is proved that the boundary conditions at the exits \eqref{e:Hughes1DBoundary-a} are mere open-end conditions.
As a consequence, the Hughes model \eqref{e:Hughes1DBoundary}, \eqref{e:Hughes1D_reformulated} is equivalent to \eqref{e:Hughes1D_reformulated} alone enforced in the whole space $\R$, but taking initial data $\bar{\rho}$ with compact support in $\mathsf{C}$ and restricting the resulting solution to $\mathsf{C}$, see \cite[Proposition~8]{AndRosSti}.
This result allows to omit the boundary conditions \eqref{e:Hughes1DBoundary} and to consider \eqref{e:Hughes1D_reformulated} in the whole space $\R$, see \cite[Definition~7]{AndRosSti}.
This leads to the following definition of entropy solution, which simplifies those introduced in \cite{DiFrancescoFagioliRosiniRussoKRM, ElKhatibGoatinRosini}.
We first need to introduce the following notation:
\[\mathcal{F}(t,x,\rho,k,\xi) \doteq \operatorname{sign}(\rho-k) \bigl( F(t,x,\rho,\xi) - F(t,x,k,\xi) \bigr),\]
where $F$ is defined as in \eqref{e:F}.

\begin{definition}
Consider a measurable initial datum $\bar{\rho} \colon \mathsf{C} \to [0, \rho_{\max}]$.
A couple $(\rho,\xi) \in \Lloc1([0,+\infty) \times \R;[0, \rho_{\max}]) \times \mathbf{Lip}([0,+\infty); \mathsf{C})$ is an entropy solution of the initial-value problem \eqref{e:Hughes1D_reformulated} if it satisfies \eqref{e:Hughes1D_reformulated-b} for a.e.\ $t>0$ as well as the entropy inequality
\begin{subequations}\label{e:entrcon}
\begin{align}\label{e:entrconA}
\iint_{\R^+\times\R} \bigl( |\rho-k| \varphi_t + \mathcal{F}(t,x,\rho,k) \varphi_x \bigr) \,\d x \,\d t
\\\label{e:entrconB}
+2 \int_{\R^+} f(k) \varphi\bigl(t,\xi(t)\bigr) \,\d t &\geq 0
\end{align}
\end{subequations}
for all $k \in [0, \rho_{\max}]$ and test functions $\varphi \in \Cc \infty\left((0,+\infty) \times \R;[0,+\infty)\right)$.
Furthermore, upon choosing a suitable representative of $\rho$, we have that $\rho$ belongs to $\C0\left([0,+\infty);\L1(\R;[0,\rho_{\max}])\right)$ with the initial condition taken in the sense $\rho(0,\cdot\,) \equiv \bar{\rho}$.

Finally, if in addition there holds $\rho \in \L \infty\bigl([0,T];\BV(\mathsf{C};[0,\rho_{\max}])\bigr)$, for all $T>0$, then we say that $(\rho,\xi)$ is a $\BV$-regular entropy solution.
\end{definition}

Condition \eqref{e:entrcon} is a Kruzhkov-type condition.
The first line \eqref{e:entrconA} originates from the Kruzhkov entropy condition \cite{Kruzhkov}.
The last line \eqref{e:entrconB} accounts for the discontinuity of the flux along the turning curve \cite{MR1927887}.
\begin{remark}\label{rem:simpler-def}
An equivalent way to provide a definition of entropy solution is to ask: 
\begin{itemize}
    \item[(a)] that the weak formulation of \eqref{e:Hughes1D_reformulated-a} hold true (this is the standard conservativity condition, implicitly contained in \eqref{e:entrcon});
    \item[(b)] that the entropy admissibility condition \eqref{e:entrcon} holds true solely for test functions vanishing at the turning curve $x=\xi(t)$, thus suppressing the line \eqref{e:entrconB}.
\end{itemize}
Therefore, the mere Rankine-Hugoniot condition is required at $x=\xi(t)$; no specific admissibility restriction is needed for the non-classical shocks at the turning curve.
This is the point of view adopted for instance in \cite{AmadoriDiFrancesco} and in \cite{AndrGirard-preprint}.
\end{remark}

\section{A Riemann-like initial datum}
\label{sec:RP}

In this section, we construct the solution to \eqref{e:Hughes1D_reformulated} with piecewise constant initial datum
\begin{equation}
\label{e:ini}
\bar\rho(x) = \left\{\begin{array}{@{}l@{\quad\hbox{if }}l@{}}
\rho_L&-1<x<0,\\
\rho_R&0\le x<1,
\end{array}\right.
\end{equation}
together with the boundary conditions \eqref{NEW:BC_Sect2.2}.
This problem was addressed in \cite[Section~4.2]{AmadoriDiFrancesco} and in \cite[Section~3.1]{ElKhatibGoatinRosini} under slightly different assumptions.

As mentioned in Subsection 2.2, we adopt the equivalent point of view of the Cauchy problem for  \eqref{e:Hughes1D_reformulated} with initial datum $\bar \rho$ extended to zero outside $\mathsf{C}$. 
The points of discontinuity of the extended initial datum are located at $x=\pm1$, $x=0$ and $x=\xi(0)$, with the last two possibly coinciding. A crucial point is a local analysis around $x=\xi(0)$ at time $t=0$, due to the change of flux along $x=\xi(t)$ and the fact that the slope $\dot\xi$ is an unknown of the problem.

\smallskip
For simplicity, let $f(\rho) \doteq \rho v(\rho) \in \C2([0,1])$ be strictly concave. 
We start by determining the initial position of the turning curve, from \eqref{e:Hughes1D_reformulated-b}, 
\eqref{e:ini}:
\[\xi(0) = -\frac{c(\rho_L)-c(\rho_R)}{2c(\rho_L)}\,.\]

If $\rho_L = \rho_R$, then the solution $\rho(t,x)$ is even in $x$ and $\xi(t)=0$ for all $t$ (see \cite[Theorem~2]{AmadoriGoatinRosini}).  
In this case, $\rho(t,0^\pm) = 0$ for all $t>0$ and two shocks are issued at $(0,0)$, on the left between $\rho_L$ and $\rho=0$ and symmetrically between $\rho=0$ and $\rho_R=\rho_L$ on the right.

Next, assume $\rho_L > \rho_R$ and note that $\xi(0) \in (-1/2,0]$.

The solution is then obtained by glueing together the solutions to the Riemann problems at $x_0\in \{0,\xi(0)\}$ and the local solutions at the boundaries $x_0\in \{-1,1\}$. Let's distinguish two cases.

\begin{itemize}
\item 
$x_0\in\{-1, 0, 1\}$: 
The sign in the flux \eqref{e:F} is respectively $(-1)$ for $x_0=-1$ and $(+1)$ for $x_0=0$, $x_0=1$.
Therefore, the Riemann problems at $x=-1$, $x=0$ and $x=1$ are standard and in each case the solution consists of a rarefaction.
\item 
$x_0=\xi(0)$: 
Here the problem rewrites as
\begin{align*}
&\rho(0,x) = \rho_L,&
&\rho_t - f(\rho)_x = 0,~  x<\xi(0),&
&\rho_t + f(\rho)_x = 0,~ x>\xi(0).
\end{align*}
Its solution may involve a non-classical shock, that is, a discontinuity that does not satisfy the Lax entropy inequalities, see \cite{MR1927887}.
\end{itemize}

Assuming that the solution $\rho$ is locally self-similar in a neighbourhood of $(t,x)=(0,\xi(0))$, then the unknowns of the problem are the constant speed $\dot{\xi}$ and the two values $\rho_\xi^\pm \dot = \rho(t,\xi(t)^\pm)$.
By \cite[Proposition~2.4]{ElKhatibGoatinRosini} the characteristic speeds enter the turning curve $x=\xi(t)$ on the side of higher density, namely
\begin{align*}
&\rho_\xi^-<\rho_\xi^+ \Longrightarrow f'(\rho_\xi^+) \leq \dot{\xi},&
&\rho_\xi^->\rho_\xi^+ \Longrightarrow -f'(\rho_\xi^-) \geq \dot{\xi}.
\end{align*}
The Rankine-Hugoniot condition along the turning curve $x=\xi(t)$ reads as
\begin{equation}
\label{e:RH}
f(\rho_\xi^+)+f(\rho_\xi^-) = \dot{\xi} (\rho_\xi^+ - \rho_\xi^-).
\end{equation}

If one of the following conditions holds:
\begin{enumerate}
\item 
$\rho_R>\hat{\rho}$,
\item 
$\rho_L\leq\hat{\rho}$ and
\[\rho_R - \rho_{\max} < \int _{\rho_R}^{\hat{\rho}} \bigl(c(\rho)-c(\rho_R)\bigr) \,\d\rho - \int_{\rho_L}^{\hat{\rho}} \bigl(c(\rho)-c(\rho_L)\bigr) \,\d\rho < \rho_{\max} - \rho_L,\]
\item 
$\rho_R\leq\hat{\rho}<\rho_L$ and\quad
$\int_{\rho_R}^{\hat{\rho}} \bigl(c(\rho)-c(\rho_R)\bigr) \,\d\rho < \rho_{\max}-\rho_L,$
\end{enumerate}
then $\rho_\xi^+ = \rho_\xi^-$. Notice that, by \eqref{e:RH}, this implies $f(\rho_\xi^\pm)=0$ and hence 
$\rho_\xi^-=\rho_\xi^+\in \{0,\rho_{\max}\}$. Moreover one has that $\rho_\xi^\pm =\rho_{\max}$ if and only if $\rho_L=\rho_{\max}$, while in the remaining cases vacuum appears, with $\rho_\xi^\pm =0$ and two shocks are issued from $x=\xi(0)$.

In all the cases not included above, the solution $\rho$ is discontinuous along $x=\xi(t)$. 
Then the bigger of the two trace values $\rho_\xi^\pm$ has to be $\rho = \rho_L$. Let $\rho = \rho_M < \rho_L$ be the other one. 
Then, two lines of discontinuity start from $x = \xi(0)$ between the states $\rho_L$, $\rho_M$ and $\rho_L$, one being a shock and the other one a non-classical shock along the turning curve.
Finally, to determine the unique value $\rho_M$ and $\xi(t)$ it is sufficient to solve the system given by \eqref{e:Hughes1D_reformulated-b} and 
\[\xi(t) = -\frac{c(\rho_L)-c(\rho_R)}{2c(\rho_L)} + \frac{f(\rho_\xi^+)+f(\rho_\xi^-)}{\rho_\xi^+ - \rho_\xi^-} \, t.\]

\smallskip
The case of $\rho_L < \rho_R$ is analogous to the one of $\rho_L > \rho_R$ because of the spatial symmetry of the problem.

\begin{remark} 
Within the initial data \eqref{e:ini}, the discontinuity point $x=0$ does not coincide with the $\xi(0)$, unless the initial data is even.  The problem of the discontinuity located at $x=\xi(0)$ is considered in \cite[Theorem~1]{AmadoriDiFrancesco}
for $v(\rho)\dot =1-\rho$ and $c(\rho)\dot =1/v(\rho)$. Under some structural conditions on the initial data, the quantity 
\[
\Psi^* \doteq \lim_{t\to 0^+} \frac 1t \left\{\int_{\xi(0)+\delta}^1-\int_{-1}^{\xi(0)-\delta} \left[ c\bigl(\rho(t,x)\bigr) - c\bigl(\rho(0,x)\bigr)\right] \, \d x \right\}
\]
is well defined and it is exploited in the construction of $\rho(t,\cdot\,)$.
Here above, $\delta>0$ is such that $\rho(t,\cdot\,)$ is well defined in $[-1,\xi(0)-\delta) \cup (\xi(0)+\delta,1]$ for small $t>0$.
In particular, for \eqref{e:ini} with $0\leq \rho_R < \rho_L < \rho_{\max} \doteq 1$ we have that:
\begin{enumerate}
\item 
If $\Psi^* \leq -2v(\rho_L)$, then there exists a unique intermediate state $\rho_M \in [0,\rho_L)$ such that $\rho_\xi$ is given by a non-classical shock along the turning curve $\xi$ between $\rho_\xi^-=\rho_L$ and $\rho_\xi^+=\rho_M$, followed by a shock between $\rho_M$ and $\rho_L$.
\item 
If $|\Psi^*| < 2v(\rho_L)$, then $\rho_\xi$ is given by a shock between $\rho_L$ and $\rho_M=0$, followed by the turning curve $\xi$ with $\rho_\xi^\pm=0$ and then by a shock between $\rho_M=0$ and $\rho_L$.
\item 
If $\Psi^* \geq 2v(\rho_L)$, then there exists a unique intermediate state $\rho_M \in [0,\rho_L)$ such that $\rho_\xi$ is given by a shock between $\rho_L$ and $\rho_M$, followed by a non-classical shock along the turning curve $\xi$ between $\rho_\xi^-=\rho_M$ and $\rho_\xi^+=\rho_L$.
\end{enumerate}
See \cite[Section~4.2]{AmadoriDiFrancesco} for an illustrative example of the conditions on $\Psi^*$.
\end{remark}

\section{Existence results}\label{sec:existence}

The first existence results for Hughes' model \eqref{e:Hughes1D_reformulated} were obtained in 2014 by Amadori, Goatin and Rosini in \cite{AmadoriGoatinRosini} for the case $v(\rho)=1-\rho$, $c(\rho)=1/v(\rho)$ and $\rho_{\max}=1$, see \cite[Theorems~2 and~3]{AmadoriGoatinRosini}.

More in detail, \cite[Theorem~2]{AmadoriGoatinRosini} deals with the \lq symmetric case\rq, namely with initial data in the space $\mathcal{S}$ of functions $\rho \in \L\infty(\mathsf{C};[0,\rho_{\max}])$ that are even, that is $\rho(-x)=\rho(x)$ for a.e.~$x\in \mathsf{C}$.
For such initial data it is proved the following existence (and uniqueness) result.

\begin{theorem}[{\cite[Theorem~2]{AmadoriGoatinRosini}}]\label{t:existence01}
Let $v(\rho)=1-\rho$ and $c(\rho)=1/v(\rho)$.
For any initial datum $\bar{\rho}$ in $\mathcal{S}$ such that $\|\bar{\rho}\|_{\L\infty} < \rho_{\max} = 1$, there exists a unique $\BV$-regular entropy solution $(\rho,\xi)$ of  Hughes' model \eqref{e:Hughes1D_reformulated} such that $\rho(t,\cdot\,) \in \mathcal{S}$ for all $t>0$.
\end{theorem}

The first step in the proof consists in showing that $\xi\equiv0$.
This, together with the Rankine-Hugoniot condition~\eqref{e:RH}, implies that $f(\rho(t,0^+))+f(\rho(t,0^-))=0$.
Hence, by the assumption $\|\bar{\rho}\|_{\L\infty} < \rho_{\max} = 1$ and the maximum principle proved in \cite[Proposition~2.5]{ElKhatibGoatinRosini}, we have $\rho(t,0^\pm)=0$.
As a final step, it is sufficient to show that the unique solution of \eqref{e:Hughes1D_reformulated} in $\mathcal{S}$ coincides on $(0,1)$ with the classical solution of the Cauchy problem for a conservation law
\[\left\{\begin{array}{@{}>{\displaystyle}l@{\quad}l@{\ }l@{}}
\rho_t + f(\rho)_x = 0, & t>0,&x\in\R,\\
\rho(0,x)=\bar{\rho}(x), &&x\in\R,
\end{array}\right.\]
where, with a slight abuse of notation, we denoted by $\bar{\rho}$ the extension of $\bar{\rho}|_{(0,1)}$ to the whole $\R$ by the value zero on $\R\setminus(0,1)$.

We recall now the existence result proposed in \cite[Theorem~3]{AmadoriGoatinRosini}, which applies to more general initial data.
Let $[x]_+ = \max\{x,0\}$, $x\in\R$.

\begin{theorem}[{\cite[Theorem~3]{AmadoriGoatinRosini}}]\label{t:WFT}
Let $v(\rho)=1-\rho$ and $c(\rho)=1/v(\rho)$.
If the initial datum $\bar{\rho}$ is in $\BV(\mathsf{C};[0,\rho_{\max}])$, is such that $\|\bar{\rho}\|_{\L\infty} < \rho_{\max} = 1$ and satisfies
\begin{equation}
\label{e:smallBV+traces}
3\|\bar{\rho}\|_{\L\infty} + \tv(c(\bar{\rho})) + [c(\bar{\rho}(-1^+))-c(1/2)]_+ + [c(\bar{\rho}(1^-))-c(1/2)]_+ < 2,
\end{equation}
then there exists a $\BV$-regular entropy solution of \eqref{e:Hughes1D_reformulated} defined globally in time.
\end{theorem}

The proof is based on the Wave-Front Tracking (WFT) algorithm \cite{DafermosWFT}, the maximum principle proved in \cite[Proposition~2.5]{ElKhatibGoatinRosini}, a convenient choice of the wave speeds of approximate rarefaction fans and the condition proposed in \cite[Theorem~1]{AmadoriDiFrancesco} to construct the solution locally at the turning point position $x=\xi(t)$.
We recall that the WFT algorithm for \eqref{e:Hughes1D_reformulated} was first analysed in \cite{MR3055243}, but only for numerical purposes.
The main difficulty in this approach is that new fronts may arise at the turning curve not only if a wave-front interacts with the turning curve, but also if two wave-fronts interact away from the turning curve.
As a result, the total variation of the solution may generically increase.
Condition \eqref{e:smallBV+traces} plays a key role as it ensures that these situations do not occur.
Details for the construction of the approximate solution via the WFT algorithm are deferred to Section~\ref{s:WFT}.

Note in passing that uniqueness results are scarce and very partial for Hughes' model. In a setting slightly more general than the one of Theorem~\ref{t:WFT} (more precisely, it is assumed that the solution is $\BV$-regular and its one-sided traces $\rho(t,\xi(t)^\pm)$ at the turning curve location are zero), uniqueness is justified in \cite[Theorem\,4]{AndRosSti} via a cumbersome Gronwall-kind argument.

Three years later, in 2017, a new approach to provide existence results was proposed by Di Francesco, Fagioli, Rosini and Russo in \cite{DiFrancescoFagioliRosiniRussoKRM}.
There, the authors pointed out that Hughes' model can been seen as two first order Lighthill-Whitham-Richards (LWR) models \cite{LWR1, LWR2} for vehicular traffic, suitably coupled at the coupling point $x=\xi(t)$, which is an inner interface splitting the whole interval $(-1,1)$ into two subintervals.
This is obvious in the case $c\equiv1$, which corresponds to pedestrians moving toward the closest exit regardless of the overall distribution (a typical behaviour in case of panic) and to the LWR model with negative velocity on $(-1,0)$ and positive velocity on $(0,1)$, see \cite[Example~1.1]{DiFrancescoFagioliRosiniRussoKRM}.
The idea was then to apply to Hughes' model the many particle approach, which was proved in 2015 to approximate the LWR model first in \cite{DiFrancescoRosini}, see also \cite{ColomboRossi, DiFrancescoFagioliRosini-BUMI, DiFrancescoFagioliRosiniRusso, HoldenRisebro02, HoldenRisebro01}.

As a first application, the authors give a lighter proof of Theorem~\ref{t:existence01} for the symmetric case, see \cite[Theorem~1.2]{DiFrancescoFagioliRosiniRussoKRM}.
Then, the authors propose in \cite[Theorem~1.3]{DiFrancescoFagioliRosiniRussoKRM} an existence result analogous to that in Theorem~\ref{t:WFT} with a condition analogous to \eqref{e:smallBV+traces}, but without involving the traces at the exits of the initial datum.
Furthermore, their results apply for more general cost functions $c$ and speed map $v$, which are assumed to satisfy {\bf (H1)} and {\bf (H2)}, respectively.
More specifically, the authors proved the following existence results.

\begin{theorem}[{\cite[Theorem~1.2]{DiFrancescoFagioliRosiniRussoKRM}}]\label{t:mpa1}
Assume that $c$ and $v$ satisfy 
{\bf (H1)} and {\bf (H2)}, respectively.
For any initial datum $\bar{\rho}$ in $\mathcal{S}$, there exists a unique $\BV$-regular entropy solution $(\rho,\xi)$ of Hughes' model \eqref{e:Hughes1D_reformulated} such that $\rho(t,\cdot\,) \in \mathcal{S}$ for all $t>0$.
\end{theorem}

\begin{theorem}[{\cite[Theorem~1.3]{DiFrancescoFagioliRosiniRussoKRM}}]\label{t:mpa2}
Assume that $c$ and $v$ satisfy 
{\bf (H1)} and {\bf (H2)}, respectively, and that $c''(\rho)>0$ for all $\rho\in[0,\rho_{\max}]$.
If the initial datum $\bar{\rho}$ is in $\BV(\mathsf{C};[0,\rho_{\max}])$, is such that $\|\bar{\rho}\|_{\L\infty} < \rho_{\max}$ and
\begin{equation}
\label{e:smallBV}
\frac{v_{\max}}{2} \left( L \, \tv(\bar{\rho}) + 3 C \right) < v(\|\bar{\rho}\|_{\L\infty}),
\end{equation}
with
\begin{align}
C &\doteq \max\bigl\{c'(\rho) \, \rho : \rho \in [0,\|\bar{\rho}\|_{\L\infty}] \bigr\},&
L &\doteq \max\bigl\{c''(\rho) \, \rho : \rho \in [0,\|\bar{\rho}\|_{\L\infty}] \bigr\},
\label{e:CL}
\end{align}
then there exists a unique $\BV$-regular entropy solution $(\rho,\xi)$ of Hughes' model \eqref{e:Hughes1D_reformulated} defined globally in time.
\end{theorem}

The details for the construction of the approximate solution via the many particle approach is deferred to Section~\ref{s:mpa}.
Here we underline the different role of the assumption $\|\bar{\rho}\|_{\L\infty} < \rho_{\max}$ in Theorem~\ref{t:WFT} and Theorem~\ref{t:mpa2}: in the former case it is required because the cost function under consideration is $c(\rho) = 1/v(\rho)$ which is not well defined at $\rho=\rho_{\max}$, whereas in the latter case it is essential to have the right-hand-side in the inequality \eqref{e:smallBV} strictly positive.

More recently, in 2021 new existence results were proposed by Andreianov, Rosini and Stivaletta in \cite{AndRosSti} for the case of a linear cost function
\begin{equation}
\label{e:linearcost}
c(\rho) = 1+\alpha\, \rho,
\end{equation}
where $\alpha\geq0$ is a parameter of the model.
The motivation for \eqref{e:linearcost} stems from the physical meaning of $\alpha$.
Indeed it corresponds to different crowd behaviours and encodes the importance given to avoid regions with high number of pedestrians.
For instance, $\alpha=0$ corresponds to panic behaviour, when people simply move towards the closest exits without avoiding crowded regions.
On the other hand, as $\alpha$ grows, so does the importance of avoiding exits chosen by a high number of pedestrians.

The approximate solution is constructed by applying a many particle approach similar to that in \cite{DiFrancescoFagioliRosiniRussoKRM}, but with two main differences, see Section~\ref{s:linear}.
First, they changed the very definition of the approximating turning curve by substituting \eqref{e:xiDPA} with \eqref{e:turning} given below.
This choice allows to link directions switching of the particles to the instants when exactly one of the particles leaves the domain $\mathsf{C}$.
This is crucial to prove rigorously the global in time existence of a discrete solution and the boundness of the evacuation time, see \cite[Theorem~20]{AndRosSti}, whereas in \cite{AmadoriGoatinRosini, DiFrancescoFagioliRosiniRussoKRM} these are implicitly assumed.
Second, unlike \cite{DiFrancescoFagioliRosiniRussoKRM}, in \cite{AndRosSti} the authors exploited the regularizing effect of the discrete version of the Oleinik's condition rather than the $\BV$-contraction property, both proved in \cite{DiFrancescoRosini}, see also \cite{DiFrancescoFagioliRosini-BUMI}.
However, to do so they need conditions on the velocity that are slightly more restrictive than {\bf (H2)} and read as follows:
\begin{enumerate}[label={\bf (H2')},leftmargin=*]
\item[{\bf (H2')}]
The speed map $v \colon [0,\rho_{\max}] \to [0,v_{\max}]$ is $\C2$, strictly decreasing, with $v(0)=v_{\max}>0$ and $v(\rho_{\max})=0$; moreover $v'(\rho)+\rho v''(\rho) \leq 0$ for all $\rho \in [0,\rho_{\max}]$.
\end{enumerate}
Note that the above condition is slightly more restrictive than requiring $f(\rho) \doteq \rho v(\rho)$ to be strictly concave. 

Their two main existence results are given in \cite[Theorems~5 and~6]{AndRosSti}.
The main novelty of these theorems is that they take into account the possible arising of non-classical shocks along the turning curve, namely discontinuities that do not satisfy the Lax entropy inequalities, see \cite{MR1927887}.
Let us stress that none of the existence results obtained in \cite{AmadoriGoatinRosini, DiFrancescoFagioliRosiniRussoKRM} considers non-classical shocks.
In fact, the assumptions on the initial data considered in Theorems~\ref{t:existence01}, \ref{t:WFT}, \ref{t:mpa1} and \ref{t:mpa2} are meant to exclude the appearance of non-classical shocks.
However, one of the main analytical features of  Hughes' model is the possible development of non-classical shocks in the solution.
Indeed, these have a physical counterpart, modelling pedestrians that switch direction during the evacuation.
In fact, pedestrians choose their direction of motion taking into account the distance from the two exits as well as avoiding crowded regions. 
As a result, due to this latter aspect, if during the evacuation pedestrians observe an increase of the crowd in front of their chosen exit as well as a decrease of the crowd at the opposite exit, then they may decide to change direction.

Both in \cite{AmadoriGoatinRosini} and \cite{DiFrancescoFagioliRosiniRussoKRM} the presence of non-classical shocks is prevented by requiring sufficient conditions, which result in considering initial data with sufficiently small total variation and imposing $\|\bar{\rho}\|_{\L\infty(\R)} < \rho_{\max}$.
On the contrary, in \cite{AndRosSti} the initial datum has arbitrarily (possibly even infinite) total variation and it can attain the maximal density $\rho_{\max}$.

We start by recalling the conditional existence result (an \lq IF-theorem\rq) under the assumption of global variation control.
We stress that in practice, this delicate assumption seems to hold for all \lq typical\rq\ choices of initial data, see the numerical tests provided in this chapter and the tests presented in~\cite{MR3149318, MR3177723, MR3619091, MR3698447, MR3055243, HUANG2009127, MR3451862, MR3277564}.

\begin{theorem}[{\cite[Theorem~5]{AndRosSti}}]\label{t:AndRosSti1}
Consider the cost function \eqref{e:linearcost}.
Assume that $v$ satisfying {\bf (H2)} is $\C2$ and such that $f(\rho)=\rho v(\rho)$ is strictly concave in $[0,\rho_{\max}]$.
Let $\bar{\rho}$ be a measurable initial datum in $\BV(\mathsf{C};[0,\rho_{\max}])$ and let $\{( \rho^{n},\zeta^{n})\}_n$ be the sequence of approximate solutions constructed in Section~\ref{s:linear}. 
Assume that for all $T>0$ there exists a constant $\hbox{\bf TV}=\hbox{\bf TV}(T)>0$ such that, for any $t\in [0,T]$ and $n\in{\mathbb{N}}$, we have
\[\tv\left( \rho^{n}(t,\cdot\,)\right) \leq \hbox{\bf TV}.\]
Then for all $T>0$, the sequence $\{( \rho^{n},\zeta^{n})\}_{n}$ converges, up to a subsequence, in $\L1((0,T)\times \mathsf{C})\times \C0([0,T])$ to a $\BV$-regular entropy solution $(\rho,\xi)$ of Hughes' model \eqref{e:Hughes1D_reformulated} defined globally in time.
\end{theorem}

We observe that the functional defined in \cite[(9)]{DiFrancescoFagioliRosiniRussoKRM} becomes trivial in the case of a linear cost \eqref{e:linearcost} and, consequently, it becomes useless.
As a result, the proof of the above theorem is quite technical, see \cite[Section~6]{AndRosSti}.

The second existence result is based on the same construction of the approximate solution described in Section~\ref{s:linear}, but it exploits a $\BV_{\rm loc}$ compactness argument via a local reduction to microscopic approximation of the LWR model, see \cite[Section~7]{AndRosSti}.

\begin{theorem}[{\cite[Theorem~6]{AndRosSti}}]\label{t:AndRosSti2}
Consider the cost function \eqref{e:linearcost}.
Assume that $v$ satisfies {\bf (H2')}.
Let $\bar{\rho}$ be a measurable initial datum and let $\{( \rho^{n},\zeta^{n})\}_n$ be the sequence of approximate solutions constructed in Section~\ref{s:linear}.
Then for all $T>0$ the sequence $\{( \rho^{n},\zeta^{n})\}_{n}$ converges, up to a subsequence, in $\L1((0,T)\times \mathsf{C})\times \C0([0,T])$ to an entropy solution $(\rho,\xi)$ of Hughes' model \eqref{e:Hughes1D_reformulated} defined globally in time.
\end{theorem}

We also recall that, as byproduct of the sharply formulated many-particle approximation scheme described in Section~\ref{s:linear}, the authors furnish two further unconditional existence results in \cite[Corollaries~33 and~34]{AndRosSti} , both excluding non-classical shocks.
The latter deals with the symmetric case and is analogous to Theorem~\ref{t:mpa1}, so it can be seen as an alternative proof of it.
The former appears to be new as it deals with initial data well separated from the origin.
More precisely, let $\mathcal{V}$ be the space of measurable functions $\rho$ in $\L \infty(\mathsf{C};[0, \rho_{\max}])$ such that $\|\rho\|_{\L1} < 2/\alpha$ and with support in $[-1,1] \setminus [-\frac{\alpha}{2}\,\|\rho\|_{\L1},\frac{\alpha}{2}\,\|\rho\|_{\L1}]$.
We have then the following existence result.

\begin{theorem}[{\cite[Corollary~33]{AndRosSti}}]\label{t:AndRosSti3}
Consider the cost function \eqref{e:linearcost}.
Assume that $v$ satisfying {\bf (H2)} is $\C2$ and such that $f(\rho)=\rho v(\rho)$ is strictly concave in $[0,\rho_{\max}]$.
Let $\bar{\rho}$ be an initial datum in $\mathcal{V}$ and let $\{( \rho^{n},\zeta^{n})\}_n$ be the sequence of approximate solutions constructed in Section~\ref{s:linear}.
Then for all $T>0$ the sequence $\{( \rho^{n},\zeta^{n})\}_n$ converges in $\L1((0,T)\times\mathsf{C})\times \C0([0,T])$ to the unique $\mathbf{BV}$-regular entropy solution $(\rho,\xi)$ of Hughes' model \eqref{e:Hughes1D_reformulated} defined globally in time and $\rho(t,\cdot\,) \in\mathcal{V}$ for all $t>0$.
\end{theorem}

\section{The wave-front tracking approach}
\label{s:WFT}

In this section we recall the construction of the approximate solution via the Wave-Front Tracking (WFT) algorithm used in \cite{ElKhatibGoatinRosini} to prove Theorem~\ref{t:WFT}.
Let $\rho_{\max} \doteq 1$, $v(\rho) \doteq 1-\rho$, $c(\rho) \doteq 1/v(\rho)$ and $\hat{\rho} \doteq 1/2$.
Fix $n \in \N$ and let $\epsilon \doteq 2^{-n}>0$.
Introduce the grid $\mathcal{G}^n \doteq \left\{i \,\epsilon : i \in \{ 0, \ldots, \epsilon^{-1}\} \right\}$ and
consider the piecewise linear function $f^n$ that interpolates linearly the points $(\rho_i, f(\rho_i))$, $\rho_i\in \mathcal{G}^n$.
Let $\bar\rho^n\in\BV(\mathsf{C};\mathcal{G}^n)$ be a piecewise constant function such that
\begin{align*}
 &|\bar\rho^n(\pm1^\mp) - \bar\rho(\pm1^\mp)| \le \epsilon, &
 &\tv\left(c(\bar\rho^n)\right) \le \tv\left(c(\bar\rho)\right) + C_1 \epsilon,
 \\
 &\|\bar\rho^n\|_{\L\infty(\mathsf{C};\R)}
 \le \|\bar\rho\|_{\L\infty(\mathsf{C};\R)} < 1,&
 &\lim_{n\to+\infty}\|\bar\rho - \bar\rho^n\|_{\L1(\mathsf{C};\R)} 
 = 0,
\end{align*}
with $C_1 = c' (\|\bar\rho\|_{\L\infty(\mathsf{C};\R)})$. 
Define $\bar\xi^n$ as the unique solution of the equation
\begin{align}\label{e:barxiWFT}
 \int_{-1}^{\bar\xi^n} c\left(\bar\rho^n(x)\right) \d x =
 \int_{\bar\xi^n}^1 c\left(\bar\rho^n(x)\right) \d x.
\end{align}
Clearly, the above formula defines $\bar{\xi}^n \in (-1,1)$ uniquely.

Let $\mathcal{R}_c$ be the classical Riemann solver, see \cite{Bressanbook}. 
Introduce the simplified Riemann solver $\mathcal{R}_s$, that replaces any rarefaction wave given by $\mathcal{R}_c$ with a rarefaction front as described below, see \eqref{e:rarefactionsspeed}.
Apply then $\mathcal{R}_s$ to solve each Riemann problem associated to the boundary $\{-1,1\}$ and to the jumps of discontinuity of $\bar\rho^n$ away from $x=\bar\xi^n$. 
Denote by $\rho_L^n$ and $\rho_R^n$ the juxtapositions of the piecewise constant functions obtained by solving with $\mathcal{R}_s$ the Riemann problems on the left of $x=\bar\xi^n$ and on the right of $x=\bar\xi^n$, respectively.
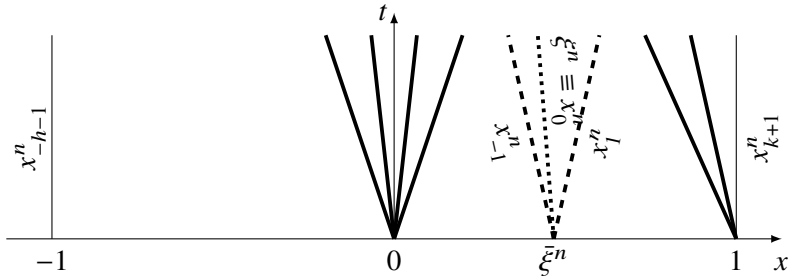
\begin{figure}[!htb]\centering
\begin{tikzpicture}[every node/.style={anchor=south west,inner sep=2pt},x=3mm, y=3mm]
\draw[-latex] (-17,0) -- (17,0) node[below] {\strut $x$};
\draw[-latex] (0,0) -- (0,10) node[left] {\strut $t$};
\draw (-15,0) node[below] {\strut $-1$} -- node[sloped, above] {\strut $x^n_{-h-1}$} (-15,9);
\draw (15,0) node[below] {\strut $1$} -- node[sloped, below] {\strut $x^n_{k+1}$} (15,9);
\foreach \a in {-3,-1,1,3}
\draw[line width=0.5mm] (0,0) -- (\a,9);
\node[below] at (0,0) {\strut $0$};
\node[below] at (7,0) {\strut $\bar\xi^n$};
\draw[line width=0.5mm, dashed] (7,0) -- node[sloped, below] {\strut $x^n_{-1}$} (5,9);
\draw[line width=0.5mm, dotted] (7,0) -- node[sloped, above left] {\strut $\xi^n\equiv x^n_{0}$} (6.3,9);
\draw[line width=0.5mm, dashed] (7,0) -- node[sloped, below] {\strut $x^n_{1}$} (9,9);
\foreach \a in {11,13}
\draw[line width=0.5mm] (15,0) -- (\a,9);

\end{tikzpicture}
\caption{Typical representation of $\rho^n$ in the case of a Riemann initial datum and obtained by juxtaposing $\rho_L^n$, $\rho_\xi^n$ and $\rho_R^n$.
Above, $\rho_L^n$ consists of the waves starting from $x=0$, $\rho_\xi^n$ consists of the waves starting from $x=\bar\xi^n$, $\rho_R^n$ consists of the waves starting from $x=1$.
The dotted line represents the turning curve $x=\xi^n(t)$.
The shock and rarefaction fronts are represented by dashed and solid thick lines, respectively.
See also \figurename~\ref{f:WFTrep}.}
\label{f:WFTnotations}
\end{figure}
One then applies \cite[Theorem~6]{AmadoriGoatinRosini}\footnote{We omit the 3~pages long \cite[Th.~6]{AmadoriGoatinRosini} to avoid overloading the Chapter with technicalities.}, which upgrades \cite[Theorem~1]{AmadoriDiFrancesco}, to construct a piecewise constant function $\rho^n_\xi$ such that if $\rho^n$ is the juxtaposition of $\rho_L^n$, $\rho_\xi^n$ and $\rho_R^n$, see \figurename~\ref{f:WFTnotations}, then the corresponding turning curve $x=\xi^n(t)$ defined by
\begin{align}\label{e:xiWFT}
 \int_{-1}^{\xi^n(t)} c\bigl(\rho^n(t,x)\bigr) \,{\d}x &=
 \int_{\xi^n(t)}^1 c\bigl(\rho^n(t,x)\bigr) \,{\d}x
\end{align}
satisfies the Rankine-Hugoniot condition \eqref{e:RH}.
As a result, the approximate solution takes the form
\begin{equation}
\label{e:approsolWFT}
\rho^n(t,x) =
\sum_{i=-h-1}^{k} \rho_{i+1/2}^n \,\mathbbm{1}_{\strut[x^n_i(t),x^n_{i+1}(t))}(x),
\end{equation}
where $\mathbbm{1}_A$ is the indicator function of $A\subset\R$, $x^n_0(t) \doteq \xi^n(t)$ is the turning curve, $x^n_i(t)$, $i \in \{ -h-1, \ldots, -1, 1, \ldots, k+1\}$, with indices $h$ and $k$  defined by imposing $-1= x^n_{-h-1}< x^n_i(t) < x^n_{i+1}(t) < x^n_{k+1} =1$, are the discontinuity lines of $\rho^n$ away from $x^n_0\equiv\xi^n$, which we call \emph{fronts}, such that
\begin{align*}
 &\rho_{i-1/2}^n \ne \rho_{i+1/2}^n \hbox{ if } i \ne0 ,&
 &\rho_{-1/2}^n = \rho_{1/2}^n\hbox{ if and only if }\rho^n_{\pm1/2} = 0 ,\\
 &\rho_{-h-1/2}^n\leq \hat{\rho} ,&
 &\rho_{k+1/2}^n\leq \hat{\rho} .
\end{align*}
When two fronts $x^n_i(t)$, $x^n_{i+1}(t)$ interact (i.e.\ $x^n_i(\bar t)=x^n_{i+1}(\bar t)$ for some $\bar t>0$ or when a front reaches the boundary), 
the approximate solution $\rho^n$ is prolonged  by applying $\mathcal{R}_s$ at interactions 
away from $x=\xi^n$ and by applying then \cite[Theorem~6]{AmadoriGoatinRosini}.
Observe that, as a result of any interaction, new fronts may originate from the turning curve, even if the interaction occurs elsewhere. 
However, the resulting approximate solution $\rho^n$ keeps the structure described above. 
Therefore, after each interaction time, we can use the same notation introduced before by rearranging the indices and by considering 
$h$ and $k$ as piecewise constant functions of time. Finally, the turning curve is prolonged by applying~\eqref{e:xiWFT} 
as long as $\rho^n$ is well defined.

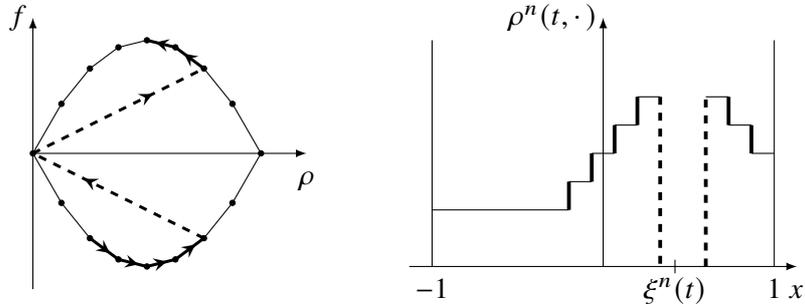
\begin{figure}[!htb]\centering
\begin{tikzpicture}[every node/.style={anchor=south west,inner sep=2pt}, x=30mm, y=60mm]

\draw[-latex] (0,0) -- (1.2,0) node[below] {\strut $\rho$};
\draw[-latex] (0,-.3) -- (0,.3) node[left] {\strut $f$};
\foreach \Point in {(0.,0.), (0.125,0.109375), (0.25,0.1875), (0.375,0.234375), (0.5,0.25), (0.625,0.234375), (0.75,0.1875), (0.875,0.109375), (1.,0.), (0.125,-0.109375), (0.25,-0.1875), (0.375,-0.234375), (0.5,-0.25), (0.625,-0.234375), (0.75,-0.1875), (0.875,-0.109375)}
    \draw[fill=black] \Point circle (1pt);

\begin{scope}[line width=0.4mm,decoration={markings, mark=at position 0.7 with {\arrow{stealth}}}]
\draw[postaction={decorate}] (0.25,-0.1875) -- (0.375,-0.234375);
\draw[postaction={decorate}] (0.375,-0.234375) -- (0.5,-0.25);
\draw[postaction={decorate}] (0.5,-0.25) -- (0.625,-0.234375);
\draw[postaction={decorate}] (0.625,-0.234375) -- (0.75,-0.1875);
\draw[postaction={decorate}, dashed] (0.75,-0.1875) -- (0,0);
\draw[postaction={decorate}, dashed] (0,0) -- (0.75,0.1875);
\draw[postaction={decorate}] (0.75,0.1875) -- (0.625,0.234375);
\draw[postaction={decorate}] (0.625,0.234375) -- (0.5,0.25);
\end{scope}

\draw (0.,0.) -- (0.125,0.109375) -- (0.25,0.1875) -- (0.375,0.234375) -- (0.5,0.25) -- (0.625,0.234375) -- (0.75,0.1875) -- (0.875,0.109375) -- (1.,0.);
\draw (0.,0.) -- (0.125,-0.109375) -- (0.25,-0.1875) -- (0.375,-0.234375) -- (0.5,-0.25) -- (0.625,-0.234375) -- (0.75,-0.1875) -- (0.875,-0.109375) -- (1.,0.);

\begin{scope}[shift={(2.5,-.25)},x=1.5mm,y=30mm]
\draw[-latex] (-17,0) -- (17,0) node[below] {\strut $x$};
\draw[-latex] (0,0) -- (0,1.1) node[left] {\strut $\rho^n(t,\cdot\,)$};
\draw (-15,0) node[below] {\strut $-1$} -- (-15,1);
\draw (15,0) node[below] {\strut $1$} -- (15,1);
\draw (-15,.25) -- (-3,.25);
\draw[line width=0.5mm] (-3,.25) -- (-3,.375);
\draw (-3,.375) -- (-1,.375);
\draw[line width=0.5mm] (-1,.375) -- (-1,.5);
\draw (-1,.5) -- (1,.5);
\draw[line width=0.5mm] (1,.5) -- (1,.625);
\draw (1,.625) -- (3,.625);
\draw[line width=0.5mm] (3,.625) -- (3,.75);
\draw (3,.75) -- (5,.75);
\draw[line width=0.5mm, dashed] (5,.75) -- (5,0);
\draw[line width=0.5mm, dashed] (9,0) -- (9,.75);
\draw (9,.75) -- (11,.75);
\draw[line width=0.5mm] (11,.75) -- (11,.625);
\draw (11,.625) -- (13,.625);
\draw[line width=0.5mm] (13,.625) -- (13,.5);
\draw (13,.5) -- (15,.5);

\draw (6.3,-.03) -- ++(0,.06);
\node[below] at (6.3,0) {\strut $\xi^n(t)$};
\end{scope}

\end{tikzpicture}
\caption{Representation of $\rho^n$ as in \figurename~\ref{f:WFTnotations} in the $(\rho,f)$-plane, left, and its profile $x\mapsto\rho^n(t,\cdot\,)$ at time $t>0$ sufficiently small, right.
The shock and rarefaction fronts are represented by dashed and solid thick lines, respectively.}
\label{f:WFTrep}
\end{figure}

To complete the construction we need to assign a travelling speed to each front.
Below, upward jumps on the left of $x=\xi^n(t)$ and downward jumps on the right of $x=\xi^n(t)$ are referred to as \emph{rarefaction fronts}, while the remaining jumps away from $x=\xi^n(t)$ are called \emph{shock fronts}, see \figurename~\ref{f:WFTrep}.
The size of the jumps is denoted by
\[\sigma_i(t) = \operatorname{sign}(i)(\rho_{i-1/2} - \rho_{i+1/2}).\]
The speeds of propagation of the turning curve and the shock fronts are obtained by imposing the Rankine-Hugoniot jump condition \eqref{e:RH}, that is
\begin{align}\nonumber
 &(\rho_{1/2} - \rho_{-1/2}) \, \dot{\xi}^n = f(\rho_{1/2}) + f(\rho_{-1/2}),\\
\nonumber
 &\dot x^n_i =\operatorname{sign}(i) \ \frac{f(\rho_{i+1/2}) - f(\rho_{i-1/2})} {\rho_{i+1/2} - \rho_{i-1/2}}
 &&\hbox{if }i\ne0 \hbox{ and }\sigma_i<0.
\intertext{Instead, any rarefaction front $x^n_i$ travels with speed}
 \label{e:rarefactionsspeed}
 &\dot x^n_i =\operatorname{sign}(i) \ \frac{q(\rho_{i+1/2}) - q(\rho_{i-1/2})} {c(\rho_{i+1/2}) - c(\rho_{i-1/2})}
 &&\hbox{if }i\ne0 \hbox{ and }\sigma_i>0,
\end{align}
where $q$ is defined by
\[q(\rho) \doteq -c(\rho) + 2 \ln\bigl(c(\rho)\bigr),\]
and is the entropy flux associated to $c$.
This choice for the speed of propagation of the rarefaction fronts allows to simplify the terms appearing in \eqref{e:xiWFT} by exploiting the fact that entropy conditions hold with an equality along any classical rarefaction.

Note that $\rho^n$ does not necessarily take values in $\mathcal{G}^n$. 
Indeed the states along the turning curve may not belong to $\mathcal{G}^n$.
However, condition \eqref{e:smallBV} ensures that this case does not occur.

\begin{remark}
The WFT scheme described below in Section~\ref{sec:sim} follows a slightly different construction. 
Indeed, there the rarefaction fronts move with the speeds prescribed by the corresponding Rankine-Hugoniot conditions. 
Also, the states along the turning curve are approximated by the closest points of the mesh $\mathcal{G}^n$. 
This choice allows to consider also non-classical shocks, see \figurename~\ref{f:04}.
\end{remark}

\section{A deterministic particles  approach}
\label{s:mpa}

In this section we recall the construction of the approximate solution via the Deterministic Particle Approximation (DPA) used in \cite{DiFrancescoFagioliRosiniRussoKRM} to prove Theorem~\ref{t:mpa2}, see also \cite{MR3644595}.

Let $\bar{\rho}$ be in $\L\infty(\mathsf{C};[0,\rho_{\max}])$. 
For a fixed $n\in \N$, set $N \doteq 2^n$ and $m \doteq 2^{-n} \, M$, where $M\doteq\|\bar{\rho}\|_{\L1}$.
Denote
\[
\bar{x}_0 \doteq \min\bigl\{\spt(\bar{\rho})\bigr\},
\]
where $\spt$ stands for the support. 
We recursively define
\begin{equation}\label{e:iniDPA}
  \bar{x}_i \doteq \inf\left\{ x > \bar{x}_{i-1} \colon \int_{\bar{x}_{i-1}}^x \bar{\rho}(y) \,\d y \geq m\right\},\quad i\in\{1,\ldots,N\}.  
\end{equation}
The above equation defines the set of $N + 1$ particles' initial positions $-1 \le \bar{x}_0 < \bar{x}_1 < \ldots < \bar{x}_{N-1} < \bar{x}_N \le 1$, with the property that the mass of the density $\bar{\rho}$ in each interval $(\bar{x}_i,\bar{x}_{i+1})$ is exactly $m$.
Introduce the local discrete initial densities
\begin{align*}
&\bar{R}_{i+1/2} \doteq \frac{m}{\bar{x}_{i+1}-\bar{x}_i},&
&i\in\{0,\ldots,N-1\},
\end{align*}
and the corresponding piecewise constant discrete initial density $\bar{\rho}^n \colon \R \to [0,\rho_{\max}]$ defined by
\[
\bar{\rho}^n(x) \doteq \sum_{i=0}^{N-1} \bar{R}_{i+1/2} \, \mathbbm{1}_{[\bar{x}_i,\bar{x}_{i+1})}(x).
\]

The initial approximated turning point $\bar{\xi}^n$ can be defined via the formula
\begin{equation}
\int_{-1}^{\bar{\xi}^n} c\bigl(\bar{\rho}^n(y)\bigr) \,\d y =
\int_{\bar{\xi}^n}^1 c\bigl(\bar{\rho}^n(y)\bigr) \,\d y.
\label{e:barxiDPA}
\end{equation}
Clearly, the above formula defines $\bar{\xi}^n \in (-1,1)$ uniquely.

By a slight modification of the initial condition, we can always assume that $\bar{\xi}^n$ does not coincide with any of the particle's initial positions.
Then, there exists $I_0\in \bigl\{0,\ldots,N\bigr\}$ such that $\bar{\xi}^n \in (\bar{x}_{I_0},\bar{x}_{I_0+1})$.
The particles on the left of $\bar{\xi}^n$ move according to a \emph{backward} follow-the-leader scheme, those on the right of $\bar{\xi}^n$ move according to a \emph{forward} one.
More precisely, we set
\begin{equation}\label{e:FtL}
\begin{cases}
\dot{x}_0(t)= -v_{\max},\\
\dot{x}_i(t)=-v\Bigl(\frac{m}{x_{i}(t)-x_{i-1}(t)}\Bigr), &i \in \{1,\ldots, I_0\},\\
\dot{x}_i(t)=v\Bigl(\frac{m}{x_{i+1}(t)-x_i(t)}\Bigr), &i \in \{I_0+1,\ldots, N-1\},\\
\dot{x}_N(t)= v_{\max},\\
x_i(0)=\bar{x}_i,&i \in \{0,\ldots, N\}.
\end{cases}
\end{equation}
We consider the corresponding local discrete densities
\begin{align*}
&R_{i+1/2}(t) \doteq \frac{m}{x_{i+1}(t)-x_i(t)},&
&i \in \{0,\ldots,N-1\}\setminus\{I_0\},
\\
&R_{i+1/2}(t) \doteq 0,&
&i \in \{-1,I_0,N\},
\end{align*}
and the corresponding piecewise constant discrete density $\rho^n \colon [0,+\infty)\times\R \to [0,\rho_{\max}]$ defined by
\begin{equation}\label{e:approsolDPA}
\rho^n(t,x) \doteq \sum_{i=0}^{N-1} R_{i+1/2}(t) \, \mathbbm{1}_{[x_i(t),x_{i+1}(t))}(x).
\end{equation}
Notice that the above density has been set to equal zero outside the particle region $[x_0(t),x_N(t))$ and around the turning point, namely in $[x_{I_0}(t),x_{I_0+1}(t))$. 
The latter in particular is simply due to a consistency with the numerical simulations, in which the computation of the turning point is made simpler in this way, see Section~\ref{sec:sim}. 
This simplifying assumption introduces a small error $m \doteq 2^{-n}M$ in the total mass.

In view of the above notation, system \eqref{e:FtL} can be written in a simpler form as follows
\begin{equation}\label{e:DevinTownsend2}
\begin{cases}
\dot{x}_i(t)=-v\bigl(R_{i-1/2}(t)\bigr),
&i \in \{0,\ldots, I_0\},
\\
\dot{x}_i(t)=v\bigl(R_{i+1/2}(t)\bigr),
&i \in \{I_0+1,\ldots, N\},\\
x_i(0)=\bar{x}_i,&i \in \{0,\ldots, N\}.
\end{cases}
\end{equation}
Notice that $R_{I_0+1/2}$ does not bias the movement of any of the particles, and this is an argument in favour of the ansatz $R_{I_0+1/2} \equiv 0$. 
The (unique) solution to the system \eqref{e:FtL} is well defined 
until the turning point does not collide with a particle.
We note that the density $R_{I_0+1/2}(t)$ is equal to zero until the turning point collides with a particle.
We shall not impose any boundary condition to the particle system \eqref{e:FtL}, and we shall follow the movement of each particle whether or not they are in $\mathsf{C} \doteq (-1,1)$. Hence, the discrete densities $R_{i+1/2}(t)$ are defined for all $t>0$ and for all $i \in \bigl\{0,\ldots,N-1\bigr\}$.

The approximate turning point $\xi^n(t)$ is implicitly uniquely defined by
\begin{equation}\label{e:xiDPA}
\int_{-1}^{\xi^n(t)} c\bigl(\rho^n(t,y)\bigr) \,\d y = \int^1_{\xi^n(t)} c\bigl(\rho^n(t,y)\bigr) \,\d y,
\end{equation}
where $\rho^n$ is the discrete density defined by \eqref{e:approsolDPA}.
Clearly $\xi^n(t)$ belongs to $\mathsf{C}$ for any $t\ge0$. 
 We emphasize that $\xi^n(0)$ does not necessarily coincide with $\bar{\xi}^n$.

We conclude this section by just highlighting the main ideas behind the proofs of Theorems~\ref{t:mpa1} and~\ref{t:mpa2}. 
Concerning Theorem~\ref{t:mpa1}, the symmetry of the initial datum implies that the discrete turning point will stuck at zero for all times, that is $\xi^n\equiv0$. 
Therefore, the particles split into two time-invariant sets, with the two particles nearest the turning point getting further and further away from each other. 
Hence, the dynamics of each group are governed by a Follow-the-Leader (FtL) model.
By the results in \cite{DiFrancescoRosini}, we then obtain convergence of $\rho^n(t,x)$ defined in \eqref{e:approsolDPA} via the FtL-Hughes particle system \eqref{e:DevinTownsend2} to the entropy solution of the Hughes model \eqref{e:Hughes1D_reformulated} as $m$ goes to zero and $N$ goes to infinity.

In proving Theorem~\ref{t:mpa2}, the first step consists in showing that by condition \eqref{e:smallBV} no particle reaches the turning curve, namely, no particle changes direction, see \cite[Proposition~2.1]{DiFrancescoFagioliRosiniRussoKRM}.
This ensures that problem \eqref{e:FtL} admits a global-in-time solution. 
In obtaining this result, a key role is played by the functional $\Upsilon(\rho) \doteq c(\rho) - c'(\rho) \, \rho$ and the following estimates
\begin{align*}
&|\dot{\xi}^n| \leq \frac{v_{\max}}{2} \Bigl(\tv\bigl(\Upsilon(\rho^n)\bigr)+3C\Bigr),&
&\dot{x}_{I_0} \leq -v(\|\bar{\rho}\|_{\L\infty}),&
&\dot{x}_{I_0+1} \geq v(\|\bar{\rho}\|_{\L\infty}),&
\end{align*}
where $C$ is defined in \eqref{e:CL}.
Notice that if $L$ is defined as in \eqref{e:CL}, then it is the Lipschitz constant of $\Upsilon$, hence $\tv\bigl(\Upsilon(\rho^n)\bigr) \leq L \, \tv(\rho^n) \leq L \, \tv(\bar{\rho})$ by the contraction estimate proven in \cite[Proposition~5]{DiFrancescoRosini}.

The second step consists in proving that $\xi^n$ converges, up to a subsequence, strongly in $\C0([0,T];\R)$ for all $T\geq0$ to some $\xi\in\C0([0,T];\R)$ and the corresponding limit turning curve $\mathcal{T} \doteq \bigl\{(t,x) \in (0,+\infty) \times \R \colon x=\xi(t)\bigr\}$ is entirely contained in the open cone 
\[\mathcal{C} \doteq \Bigl\{(t,x) \in (0,+\infty) \times \R : |x-\bar\xi| < \frac{v_{\max}}{2} \bigl( L \, \tv(\bar{\rho}) + 3 C \bigr) \, t\Bigr\}.\]

Then, since no particle is placed in $\mathcal{C}$, the discrete density $\rho^n(t,x)$ defined in \eqref{e:approsolDPA} converges to zero strongly in $\Lloc1(\mathcal{C})$. 
Applying then the results in \cite{DiFrancescoRosini}, one can prove that on $[\delta,+\infty)\times [0,+\infty)$, with $\delta>0$, the discrete density $\rho^n(t,x)$ converges strongly in $\Lloc1$ towards a function $\rho_R\in \L\infty([\delta,+\infty)\times \R)$ satisfying the Kruzhkov's entropy condition \cite{Kruzhkov} for the conservation law $\rho_t+f(\rho)_x = 0$. 
Similarly, on $[\delta,+\infty)\times (-\infty,0]$, with $\delta>0$, we have that $\rho^n(t,x)$ converges strongly in $\Lloc1$ towards a function $\rho_L\in \L\infty([\delta,+\infty)\times \R)$ satisfying the Kruzhkov's entropy condition \cite{Kruzhkov} for the conservation law $\rho_t-f(\rho)_x = 0$. 
As a consequence, $\rho^n(t,x)$ converges to
\[\rho(t,x)=
\begin{cases} \rho_L(t,x) & \text{if } x<\xi(t), \\
\rho_R(t,x) & \text{if } x>\xi(t).
\end{cases}\]
At last, it is easy to check that $\rho$ defined above satisfies the initial condition, hence $(\rho,\xi)$ is the entropy solution of the Hughes model \eqref{e:Hughes1D_reformulated}.

\begin{remark}
Clearly, equations \eqref{e:barxiDPA}, \eqref{e:approsolDPA} and \eqref{e:xiDPA} introduced for the DPA scheme are similar to \eqref{e:barxiWFT}, \eqref{e:approsolWFT} and \eqref{e:xiWFT} introduced for the WFT algorithm, respectively.
However, let us underline the main difference between \eqref{e:approsolDPA} and \eqref{e:approsolWFT}: in the former equation the extremes of the sum do not depend on time, whereas in the latter equation both the extremes are in general piecewise constant functions of time.
As a result, the DPA scheme is somehow simpler than the WFT algorithm.
\end{remark}

\section{The case of a linear cost function}
\label{s:linear}

In this section we recall the Follow-the-Leader (FtL) Hughes particle model proposed in \cite{AndRosSti} to construct the approximate solution used to prove Theorems~\ref{t:AndRosSti1}, \ref{t:AndRosSti2} and~\ref{t:AndRosSti3}.
Recall that the authors consider in \cite{AndRosSti} a linear running cost function $c(\rho)$, that is
\begin{equation}\label{e:cost}
c( \rho ) = 1+\alpha \, \rho
\end{equation}	
with $\alpha\ \geq \ 0$.
Let us stress that the first advantage in choosing a linear cost function is the opportunity to reproduce different crowd behaviours with the same model, by just letting vary the value of the parameter $\alpha$. 
Furthermore, we can assign a physical meaning to $\alpha$: it measures the importance given to avoid regions with a high number of pedestrians.
In fact, taking $\alpha=0$ corresponds to a panic behaviour, when people simply move towards the closest exit.
On the other hand, as $\alpha>0$ grows so does the importance of avoiding exits attracting a high number of pedestrians.

The strategy of existence analysis proposed in \cite{AndRosSti} is similar to that already used in previous works on the many-particle approximation of the one-dimensional Hughes model \eqref{e:Hughes1D_reformulated}, see \cite{MR3644595, DiFrancescoFagioliRosiniRusso}.
The only (crucial) difference is the definition of the approximate turning curve.
As a result, the construction of the piecewise constant discrete density $\rho^n \colon [0,+\infty)\times\R \to [0,\rho_{\max}]$ proposed in \cite{AndRosSti} is analogous to that one proposed in \cite{DiFrancescoFagioliRosiniRusso} and already described in Section~\ref{s:mpa}: in  simple words, it is sufficient to replace $\xi^n$ defined by \eqref{e:xiDPA} with $\zeta^n$ defined by \eqref{e:turning} given below.
For this reason, below we will not recall the Deterministic-Particle-Approximation (DPA) used in \cite{AndRosSti} to construct an approximate solution, see \cite[Section~6]{AndRosSti}, but we rather describe the FtL Hughes particle model, see \cite[Section~5]{AndRosSti}.

Fix $M>0$, $n\in \mathbb{N}$ and $-1\leq \overline{x}_0<\overline{x}_1<\dots<\overline{x}_N\leq 1$, with $N\doteq2^n$, satisfying
\[\overline{x}_{i+1}-\overline{x}_{i} \geq \frac{m}{\rho_{\max}},\]
for all $i\in\{ 0,\dots, N-1\}$, where $m \doteq M/N$. 

The time evolution in the whole of $\mathbb{R}$ of the particle system $x_{0}(t),\ldots,x_{N}(t)$ is described by the FtL system
\begin{equation}
\label{e:FtLlinear}
\left\{\begin{array}{@{}l@{\qquad}l@{\quad}l@{}}
\dot{x}_{i}(t) = -v\left(R_{i-1/2}(t)\right) &\hbox{ if } x_{i}(t) < \zeta^{n}(t),& i \in \{ 0 ,\dots, N \},
\\
\dot{x}_{i}(t) = v\left(R_{i+1/2}(t)\right) &\hbox{ if } x_{i}(t) \geq \zeta^{n}(t),& i \in \{ 0 ,\dots, N \},
\\
x_{i}(0) = \overline{x}_{i}, && i \in \{ 0 ,\dots, N \},
\end{array}\right.
\end{equation}
where
\begin{equation}
\label{e:Plini}
R_{i+1/2}(t) \doteq 
\frac{m}{x_{i+1}(t)-x_{i}(t)},
\qquad i \in \{ -1,\dots, N \},
\end{equation}
and
\begin{align}
\label{e:IBuiltTheSky}
&x_{-1}(t)\doteq -\infty,&
&x_{N+1}(t)\doteq +\infty.
\end{align}
Notice that by \eqref{e:Plini} and \eqref{e:IBuiltTheSky} we have $R_{-1/2} \equiv 0$ and $R_{N+1/2} \equiv 0$, therefore $v(R_{-1/2}) \equiv v_{\max}$ and $v(R_{N+1/2}) \equiv v_{\max}$.
As for the ODE system \eqref{e:FtL} introduced in the previous section, the ODE system \eqref{e:FtLlinear} needs to be closed by providing the dynamics of the turning point $\zeta^{n}(t) \in \mathbb{R}$. 
In place of $\xi^n(t)$ implicitly defined by \eqref{e:xiDPA} in the previous section, here we consider $\zeta^{n}(t)$ implicitly (uniquely) determined by 
\begin{equation}
\label{e:turning}
Z_-\left(t,\zeta^{n}(t)\right) = Z_+\left(t,\zeta^{n}(t)\right),
\end{equation}
where $Z_\pm \colon [0,+\infty) \times \mathbb{R} \to \mathbb{R}$ are the piecewise linear continuous functions defined by
\begin{align*}
Z_-(t,x) &\doteq 
\begin{cases}\displaystyle x +1 + \alpha \int_{x_{I_-}(t)}^{x} \rho^{n}(t,y) \, {\rm{d}} y
&\begin{minipage}[t]{.4\linewidth}
if $\exists\, I_- \in \{ 0,\dots, N \}$ such that\\$x_{I_--1}(t) \leq -1 < x_{I_-}(t) < x$,
\end{minipage}
\\[7pt] \displaystyle
x + 1 &\hbox{otherwise},
\end{cases}
\\
Z_+(t,x) &\doteq 
\begin{cases}\displaystyle 1 - x + \alpha \int_{x}^{x_{I_+}(t)} \rho^{n}(t,y) \, {\rm{d}} y
&\begin{minipage}[t]{.4\linewidth}
if $\exists\, I_+ \in \{ 0,\dots, N \}$ such that\\$x < x_{I_+}(t) < 1 \leq x_{I_++1}(t)$,
\end{minipage}
\\[7pt] \displaystyle
1 - x &\hbox{otherwise},
\end{cases}
\end{align*}
with $\rho^{n} \colon (0, +\infty) \times \mathbb{R} \to [0, \rho_{\max}]$ being the discrete density
\begin{equation}
\label{e:disdens}
\rho^{n}(t,x) \doteq \sum_{i = 0}^{N-1}R_{i+1/2}(t) \, \mathbbm{1}_{\left[x_{i}(t),x_{i+1}(t)\right)}(x).
\end{equation}

To sum up, the many-particle approximation consists in the ODE system \eqref{e:FtLlinear}-\eqref{e:IBuiltTheSky}, which features discontinuities in the state variable $(x_0,\dots,x_N)$ driven by the function $\zeta^n$ implicitly determined by relations \eqref{e:turning}-\eqref{e:disdens}. 
The rigorous notion of solution for this microscopic model is given in the following definition.
\begin{definition}
\label{d:dsol}
We say that an $(N+2)$-tuple $\left((x_0,\dots,x_N),\zeta^n\right)$ of functions defined on $[0,\tau)$ (for some $\tau\in (0,+\infty]$) is a solution to \eqref{e:FtLlinear}-\eqref{e:disdens} if it satisfies \eqref{e:FtLlinear}-\eqref{e:disdens} and has the following regularity:
\begin{itemize}
\item [(i)] $x_i$, $i\in \{ 0,\dots, N \}$, and $\zeta^n$ are piecewise $\C1$ on $[0,\tau)$. More precisely, there exists $H_{\rm sw}\in \mathbb{N}$ and times $\{t_h\}_{h\in \{ 1,\dots, H_{\rm sw} \}}$, $t_1<t_2<\dots<t_{H_{\rm sw}}<\tau$, such that, upon setting $t_0=0$ and $t_{H_{\rm sw}+1}=\tau$, the restriction of each of these functions to the time intervals $(t_h,t_{h+1})$ can be extended on $[t_h,t_{h+1})$ as a $\C1$-function.
\item [(ii)] $x_i$, $i\in \{ 0,\dots, N \}$, are continuous on $[0,\tau)$, while their derivatives $\dot x_i$ and the function $\zeta^n$ are normalized by the left-continuity at the times $t_h$, $i\in \{ 1 ,\dots, H_{\rm sw} \}$.
\end{itemize} 
\end{definition}
The new definition of the turning curve $x=\zeta^n(t)$, which is deeply linked to the choice of the cost function \eqref{e:cost}, allows to better highlight the microscopic counterpart of the arise of non-classical shocks (see \cite{MR1927887}) for the Hughes model: pedestrians may switch direction during the evacuation of $\mathsf{C}$.
In fact, pedestrians choose their direction of motion according to a weighted distance encoding the overall distribution of the crowd in $\mathsf{C}$.
Therefore, pedestrians choose their path towards the fastest exit, taking into account both the distance from the two exits as well as avoiding crowded regions.
Moreover, the relevance given to the first or the second factor depends on the value of the parameter $\alpha$. 
This leads to a many-particle dynamic for which the instants of particles' interactions with the turning curve are sharply captured, and this allows to get a rigorous construction of the unique global in time solution to the many-particle system. 
In fact, the \emph{a priori} analysis of solutions of \eqref{e:FtLlinear}-\eqref{e:disdens} carried out in \cite[Section~5]{AndRosSti} implies that, once a particle leaves $\mathsf{C}$, it cannot re-enter, that is it remains outside $\mathsf{C}$; moreover, at most just one particle can interact with the turning curve, and this can occur only if exactly at the same instant of time exactly one particle leaves $\mathsf{C}$.
This is proved by showing that the turning curve $x=\zeta^n(t)$ has a discontinuity only in the case exactly one particle leaves $\mathsf{C}$, moreover the turning curve can cross a particle trajectory only by jumping across it and, in this case, it crosses exactly one particle trajectory. 
This allows to link direction switches of particles (i.e., crossings of particles' paths with the turning curve) to the instants when exactly one of the particles leaves the domain $\mathsf{C}$. 
The combination between these peculiar features with the usual (in the context of the FtL approximation) discrete maximum principle culminates in the following result.

\begin{theorem}[{\cite[Theorem~20]{AndRosSti}}]
System \eqref{e:FtLlinear}-\eqref{e:IBuiltTheSky} coupled to \eqref{e:turning}-\eqref{e:disdens} admits a unique global solution in the sense of Definition \ref{d:dsol}.
Furthermore, there exists $T>0$ such that $\rho^n(t,\cdot\,)\equiv0$ in $\mathsf{C}$ for any $t\geq T$.
\end{theorem}

By the above theorem, the microscopic evacuation time $T_{\rm mic} = \inf\{t>0 : \rho^n(t,\cdot\,)\equiv0 \hbox{ in }\mathsf{C}\}$ is bounded.
Let us emphasize that this property is not proved for the approximate solutions constructed in \cite{AmadoriGoatinRosini} or in \cite{DiFrancescoFagioliRosiniRussoKRM}.

Furthermore, in \cite[Section~8]{AndRosSti} it is studied how the parameter $\alpha$ impacts on $T_{\rm mic} = T_{\rm mic}(\alpha)$, showing that it may have infinitely many discontinuities and a global minimum, see \cite[Figure~6]{AndRosSti}.

\section{Fixed-point existence strategy}
\label{s:fixed-point}

The fixed-point approach to Hughes' problem in one space dimension is the subject of the very recent work \cite{AndrGirard-preprint}. For the original Hughes' model, this approach   yields existence under the same assumption of linear costs as \cite{AndRosSti}, and under weaker restrictions on the velocity profile $v$; however, this approach also allows to consider  more general \lq capacity drop\rq\ behaviour at exits as introduced in \cite{ADR2014} (instead of the standard exit conditions). For generalised Hughes' models where the uniform Lipschitz continuity is guaranteed for a turning curve $\xi$ computed from a given density $\rho$ (two examples of such models, involving memory and relaxation effects, are given in \cite{AndrGirard-preprint}), existence of a solution follows from the method of \cite{AndrGirard-preprint}, for general costs. Let us briefly describe these results.

Fix a finite time horizon $T>0$. We reformulate the one-dimensional Hughes' problem cast under the form \eqref{e:Hughes1D_reformulated} - and suggest formulating a wide family of its abstract generalisations - by considering a solution of \eqref{e:Hughes1D_reformulated} as a fixed point of the composition $\mathcal{S}_0\circ\mathcal{I}_0$ of two operators:
\begin{gather}\label{eq:S0-operator}
    \text{$\mathcal{S}_0$ maps a given $\xi\in \Lip([0,T])$ to $\rho$ solving \eqref{e:Hughes1D_reformulated-a}, \eqref{e:Hughes1D_reformulated-c}};
\\\label{eq:I0-operator}    
    \text{$\mathcal{I}_0$ maps a given $\rho\in \L1((0,T) \times \mathbb{R})$ to $\xi\in \C0([0,T])$ solving \eqref{e:Hughes1D_reformulated-b}}.
\end{gather}
Upon replacing the operators $\mathcal{S}_0$, respectively $\mathcal{I}_0$, by $\mathcal{S}$ (that may correspond to the solution of a variant of \eqref{e:Hughes1D_reformulated-a}, \eqref{e:Hughes1D_reformulated-c}), respectively by $\mathcal{I}$ (that may correspond to a different modelling of collective dynamics towards exits), one obtains a wide family of generalised Hughes' models tractable within this fixed-point formalism. The alternative choices for the solver $\mathcal{S}$ may correspond to different exit conditions (see \cite[Section 4]{AndrGirard-preprint}, which we outline below) and/or to a more general expression of the flux $F=F(t,x,\rho,\xi)$ than the one given in \eqref{e:F} (see \cite[Remark 1.4]{AndrGirard-preprint} for the case with directional anisotropy of agents' movement corresponding, e.g., to a slanted corridor).
The alternative choices for the solver $\mathcal{I}$ may reflect memory and relaxation effects (see \cite[Sections~3.2 and~3.3]{AndrGirard-preprint}, which we outline below).

In order that the above introduced fixed-point problem be consistent, one needs to ensure that the above operators $\mathcal{S}_0$, $\mathcal{I}_0$ are well defined (in particular, that they are single-valued) and that they map between adequately chosen functional spaces. The definitions \eqref{eq:S0-operator}, \eqref{eq:I0-operator} highlight the choice of Banach spaces $\L1((0,T) \times \mathbb{R})$, $\C0([0,T])$ and $\Lip([0,T])$, endowed with their standard norms, in the construction.
On the one hand, one needs that the discontinuous-flux conservation law \eqref{e:Hughes1D_reformulated-a} admits a unique admissible solution (see in particular Remark~\ref{rem:simpler-def}; cf.~\cite[Definition~1.1 and Theorem~2.1]{AndrGirard-preprint} for details) and that the solver $\mathcal{S}_0$ is continuous (with respect to the well-chosen topologies). These claims hold true under the Lipschitz continuity of $\xi$ and mild assumptions on $f$. On the other hand, in order to be able to compose the two operators $\mathcal{S}_0$ and $\mathcal{I}_0$ and apply fixed-point arguments, one needs that,  for an appropriately defined convex closed bounded subset $B$ of $\L1((0,T) \times \mathbb{R})$,  the function $\xi=\mathcal{I}_0[\rho]$ belongs to $\Lip([0,T])\subset \C0([0,T])$ whenever $\rho\in B$. Note that the embedding  of $\Lip([0,T])$ into $\C0([0,T])$ is compact, which permits to apply the Schauder fixed-point theorem.
The cornerstone of the fixed-point formulation of the Hughes model, that pre-determines the above choice of the functional framework and restricts the main result to the case of a linear cost \eqref{e:linearcost}, is the following couple of observations proved in \cite{AndrGirard-preprint}.
\begin{proposition}\label{prop:Theo}
The operator $\mathcal{S}_0$ defined from $\Lip([0,T])$ (endowed with the norm of the larger space $\C0([0,T])$) to  $\L1((0,T);[0,\rho_{\max}])$ is continuous.\\
The operator $\mathcal{I}_0$ is continuous from $\L1((0,T);[0,\rho_{\max}])$ to $\C0([0,T])$. Moreover, if the cost $c$ is of the form $c(\rho) \doteq 1+\alpha\rho$, $\alpha>0$, and $\rho$ verifies the property
\begin{align}\label{eqRhoSemiContinuity}
\exists C>0 :\>&\forall a,b \in \mathbb{R}, \, \forall s,t \in [0,T] \textrm{, }& 
&\left| \int_a^b \bigl(\rho(t,x) - \rho(s,x)\bigr) \, \d x \right| \leq C |t-s|,
\end{align}
then $\xi=\mathcal{I}_0[\rho]$ belongs to $\Lip([0,T])$ and the Lipschitz constant of $\xi$ does not exceed the universal bound $\alpha C$. 
\end{proposition}

This permits us to apply the following theorem (\cite[Theorem~1.9]{AndrGirard-preprint}) pertaining to the solver $\mathcal{S}_0$ in \eqref{eq:S0-operator} and to an abstract operator $\mathcal{I}$ serving to compute the turning curve $\xi$ from the density $\rho$.
\begin{theorem}\label{thMainCorollary}
  Let $\bar \rho$ be a datum supported in $\mathsf{C}$ with values in $[0,\rho_{\max}]$. 
  Let $B$ be a convex closed bounded subset of $\L1((0,T) \times \mathbb{R})$
  and 
  \[\mathcal{I} \colon  (B,\|\cdot\|_{\L1((0,T)\times \mathbb{R})}) \to(\C0([0,T]), \|\cdot\|_{\infty})\]
  be a continuous operator. Assume that $f$ is non-degenerate in the sense of \cite{MR2592291, MR1869441}, where $f$ is involved in the definition \eqref{e:F} of the flux of the conservation law \eqref{e:Hughes1D_reformulated-a}.
  If there exists $r>0$ such that:
  \begin{align*}
    &\mathcal{I}(B) \subset B_{\Lip}(0,r) \doteq \left\{\xi \in \Lip([0,T]) : \|\dot{\xi}\|_{\infty} + \|\xi\|_{\infty} \leq r \right\} ,
    \\
    &\forall \xi \in B_{\Lip}(0,r) \textrm{, the unique admissible solution to  \eqref{e:Hughes1D_reformulated-a}, \eqref{e:Hughes1D_reformulated-c} is in } B , 
  \end{align*}
  then there exists a solution $(\rho,\xi)$ to the generalised Hughes' problem \eqref{e:Hughes1D_reformulated-a}, $\xi=\mathcal{I}[\rho]$ (generalising \eqref{e:Hughes1D_reformulated-b} that corresponds to the case $\mathcal{I}=\mathcal{I}_0$), \eqref{e:Hughes1D_reformulated-c}.
\end{theorem}
To apply this result to the original Hughes' model (i.e., $\mathcal{I}=\mathcal{I}_0$ given by \eqref{eq:I0-operator}), it remains to use Proposition~\ref{prop:Theo} and observe that the set $B=B_1$, with
\[
  B_{1} \doteq \bigl\{
  \rho : \|\rho\|_{\L1((0,T)\times\mathbb{R})}\leq T\|\bar \rho\|_{\L1} \textrm{ s.t.\ } 0 \leq \rho \leq \rho_{\max} \textrm{ and } \rho \textrm{ verifies \eqref{eqRhoSemiContinuity}} \bigr\},
\]
fulfils the requirements of Theorem~\ref{thMainCorollary} (see \cite[Section 3.1]{AndrGirard-preprint}). Thus Proposition~\ref{prop:Theo} and Theorem~\ref{thMainCorollary} provide a second existence proof in the context of affine costs. We stress that this proof is non-constructive (although constructive splitting arguments can be developed instead of the fixed-point arguments). Moreover, the first proof given in \cite{AndRosSti} offers appealing microscopic foundations to the macroscopic Hughes' model. The main asset of the approach of Theorem~\ref{thMainCorollary} developed in \cite{AndrGirard-preprint} is its flexibility, while the DPA approach of \cite{AndRosSti}, as well as the WFT approach of \cite{ElKhatibGoatinRosini}, require heavy adaptations if ingredients of the model change even slightly. Let us briefly discuss applications and generalizations of Theorem~\ref{thMainCorollary}.

First, the open-end exit conditions for $\rho$, implicitly contained in the formulation \eqref{e:Hughes1D_reformulated} (see \cite[p.~220]{DiFrancescoFagioliRosiniRussoKRM} and \cite[Section~3]{AndRosSti}), can be replaced by exit behaviour of the capacity drop kind, following the ideas and techniques put forward in \cite{ADRR2016, ADR2014}. Focusing on the exit situated at $x=1$ (the case of the other exit is analogous), in addition to \eqref{e:Hughes1D_reformulated-a}, \eqref{e:Hughes1D_reformulated-c}, we consider the non-local point constraint
\begin{equation}\label{seqExitR}
     f\bigl(\rho(t, 1)\bigr) \leq g \left( \int_\sigma^1 w_1(x) \rho(t,x) \, \d x \right) 
\end{equation}
for a given $g\in \Lip([0,\rho_{\max}];[0,f(\bar \rho)])$ and a given weight $w_1\in \Lip((-\infty,1])$ with support in some compact vicinity $[\sigma,1]$ of the exit. Precise definition of admissible solution $\rho$ and the existence result (given $\xi\in\Lip([0,T])$) for this original discontinuous-flux, non-locally constrained at both exits conservation law can be found in \cite[Section~4 and Appendix]{AndrGirard-preprint}. Denoting by $\mathcal{S}_g$ the associated solver, one easily obtains the analogous of Theorem~\ref{thMainCorollary} with $\mathcal{S}_g$ replacing the basic solver $\mathcal{S}_0$. Even more complex dynamics at exits, exhibiting  self-organisation features  \cite{AndrSylla-FVCA}, can be considered in the same way to replace $\mathcal{S}_0$.

Second, fixing either $\mathcal{S}_0$ in \eqref{eq:S0-operator} or $\mathcal{S}_g$ with the additional exit constraint  \eqref{seqExitR}, we can replace the operator $\mathcal{I}_0$ resolving \eqref{e:Hughes1D_reformulated-b} by different (though closely related) operators that regularize the dynamics of $\xi$, making trivial the Lipschitz bound on $\xi=\mathcal{I}[\rho]$. We stress that in this case, the analogous of Proposition~\ref{prop:Theo} holds true without the restrictive assumption \eqref{e:cost} of linear costs. For the first example, we introduce a memory effect via a subjective density
\[
\mathcal{R}[\rho(\,\cdot,x)](t) \doteq
\delta \int_{-\infty}^t \rho(s,x) \, e^{-\delta (t-s)} \,\d s
\]
(where $\rho$ is extended by the initial datum $\bar\rho$ for all $t<0$).
Instead of \eqref{e:Hughes1D_reformulated-b}, we define the operator $\mathcal{I}_\delta \colon \rho\to \xi$ by
\[
    \int_{-1}^{\xi(t)} c\bigl(\mathcal{R}[\rho(\,\cdot,x)](t)\bigr) \,\d x  = \int_{\xi(t)}^1 c\bigl(\mathcal{R}[\rho(\,\cdot,x)](t)\bigr) \,\d x. 
\]
The operator $\mathcal{I}_\delta$ possesses the properties required in Theorem~\ref{thMainCorollary} for the straightforward choice $B=B_2 \doteq \bigl\{
  \rho : \|\rho\|_{\L1((0,T)\times\mathbb{R})}\leq T\|\bar \rho\|_{\L1} \textrm{ s.t.\ } 0 \leq \rho \leq \rho_{\max} \bigr\}$.
For the second example, we define $\xi$ via a relaxation mechanism. The simplest variant is the ODE problem 
\[\begin{cases}\displaystyle
  -\epsilon \, \dot{\xi}(t) = \int_{\xi(t)}^1 c\bigl(\rho(t,x)\bigr) \,\d x - \int_{-1}^{\xi(t)} c\bigl(\rho(t,x)\bigr) \,\d x, 
  \\\displaystyle
  \int_{\xi(0)}^1 c\bigl(\rho_0(x)\bigr) \,\d x - \int_{-1}^{\xi(0)} c\bigl(\rho_0(x)\bigr) \,\d x = 0, 
\end{cases}\]
which defines  $\tilde{\mathcal{I}}_\epsilon \colon \rho\to \xi$ consistent with the choice $B=B_2$ in the context of Theorem~\ref{thMainCorollary} (or of its analogous involving $\mathcal{S}_g$  instead  of $\mathcal{S}_0$).

\section{Simulations}
\label{sec:sim}

The analytical results on Hughes' model illustrated in the previous sections were coupled in the literature with different numerical schemes, see \cite[Section~5]{ElKhatibGoatinRosini} and \cite{MR3055243} for the Wave-Front Tracking (WFT) scheme, and \cite[Section~3]{DiFrancescoFagioliRosiniRussoKRM} and \cite[Section~8]{AndRosSti} for the Deterministic-Particle-Approximation (DPA) algorithm \cite{DiFrancescoRosini}. 
Indeed, the introduced numerical schemes can be viewed both as analytical and numerical tools.

Several contributions can be found in the literature concerning the numerical study of  Hughes' model. In \cite{HUANG2009127}, the authors introduced a WENO scheme for the scalar conservation law and a fast sweeping method for the eikonal equation.
In \cite{MR3277564}, using a mixed finite volume method, a comparison between solutions of Hughes’ model and a second-order model was presented using extensive numerical experiments, including the case of obstacles in the interior of the domain.
The study in~\cite{carrillo} addresses the case of Hughes' model with limited local vision  both in one and two dimensions. 
A semi-Lagrangian scheme was used in~\cite{MR3698447}  to solve  both the stationary Hamilton-Jacobi equation and the regularised transport equation described in Section~\ref{s:regularised} on  bounded domains.

All the aforementioned results focus on the two-dimensional model \eqref{e:Hughes2D}, \eqref{e:HughesIni}, \eqref{e:HughesBoundary}, since it represents the most interesting case from the application point of view.
However, in this section, we limit our presentation to the one dimensional case, where analytical results are available, showing some numerical tests performed using the WFT scheme and the DPA algorithm. \\
In all the reported examples, we choose the velocity and cost functions as
\begin{align*}
&v(\rho)\doteq1-\rho,&
&c(\rho)\doteq1/v(\rho).
\end{align*}
We show the time evolution of the discrete densities constructed through the two approximations in the domain $\mathsf{D}=\left(-1,1\right)$. 
We briefly present the methods below.

\medskip

\begin{figure}[!ht]\centering
\includegraphics[scale=0.357]{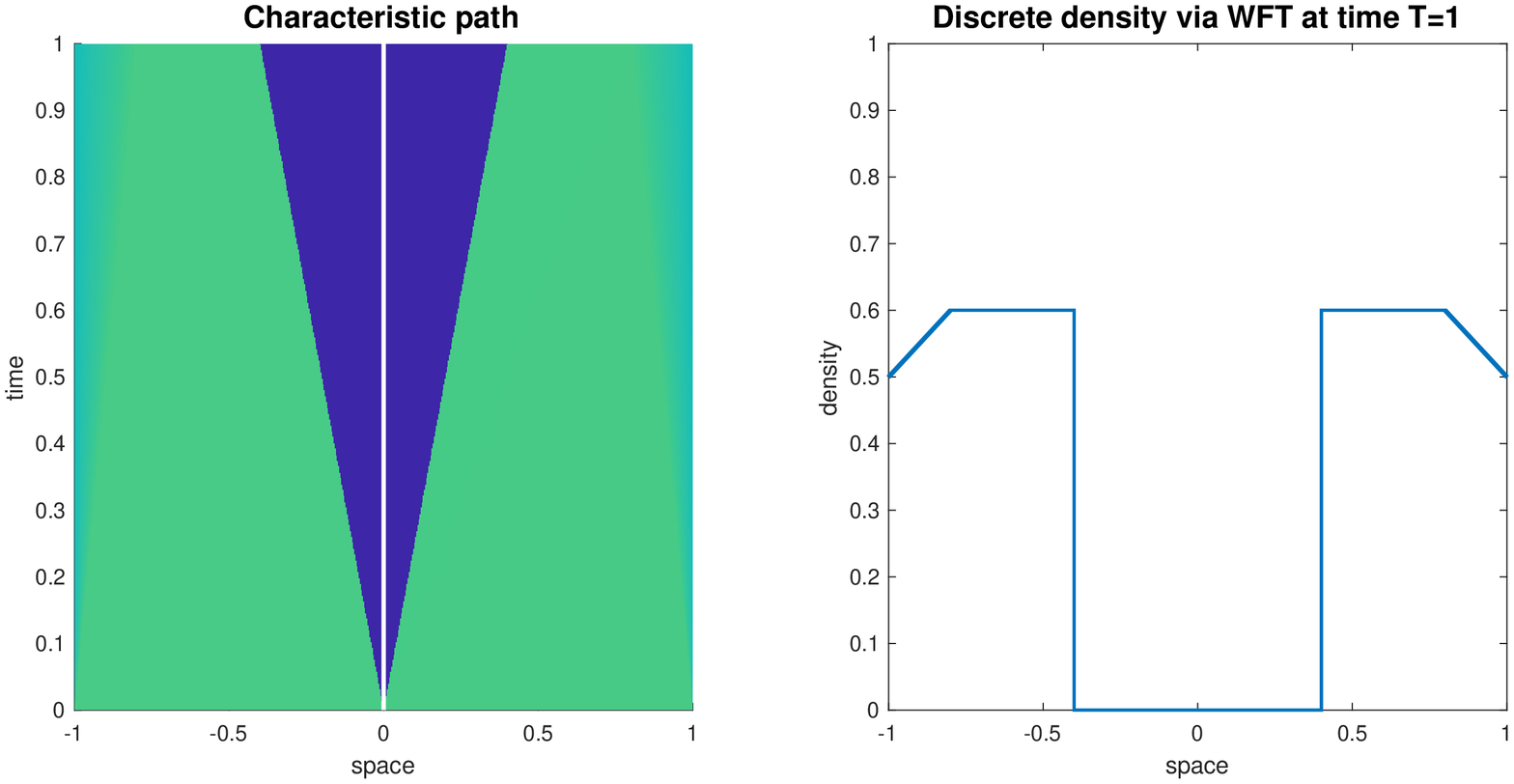}\\
\includegraphics[scale=0.357]{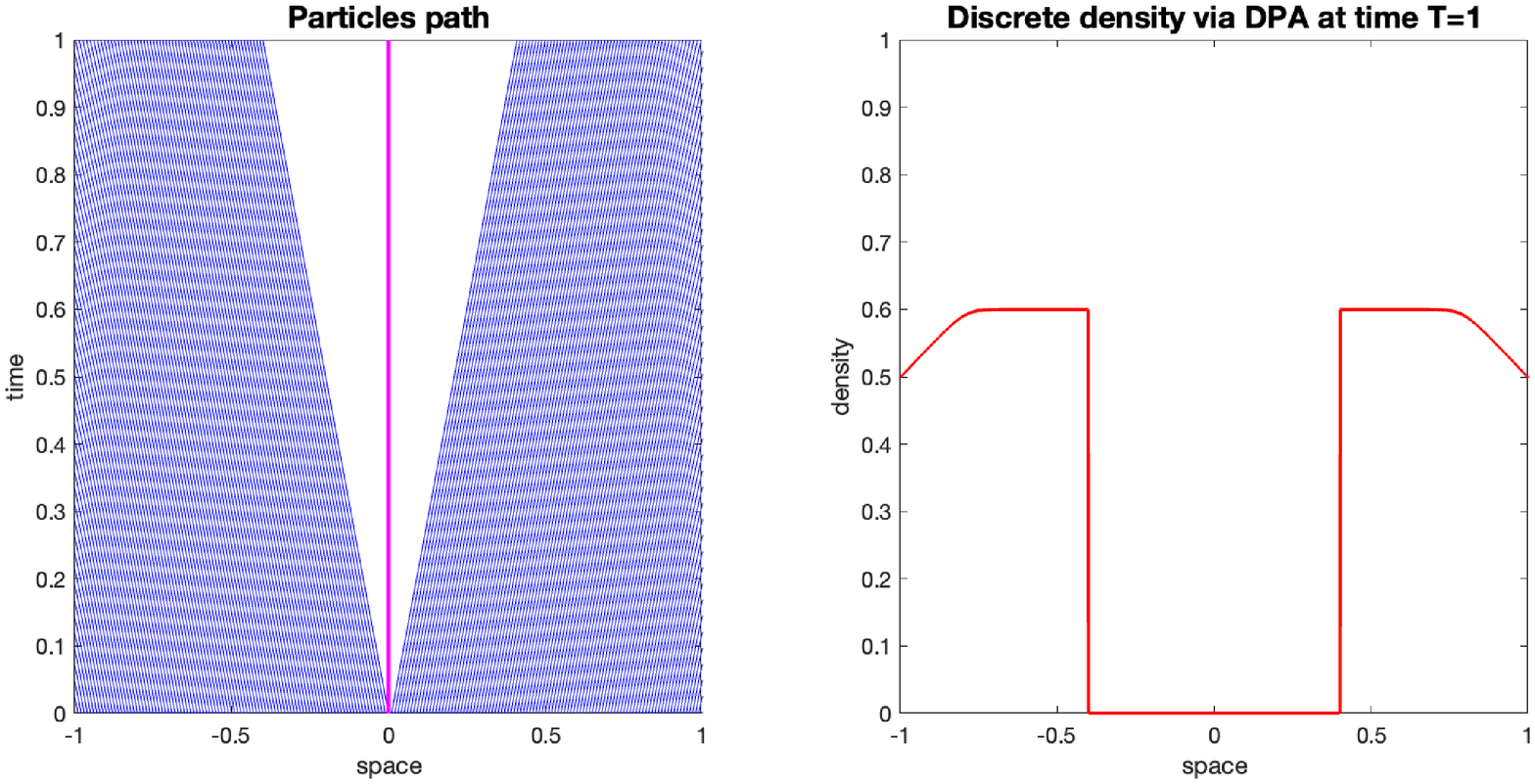}
\caption{
Approximate solutions of Hughes' model with even initial  condition~\eqref{e:infig1} constructed by the WFT scheme (top) and DPA method (bottom). The spatio-temporal evolution is depicted on the left, the corresponding density profile at time $t=1$ on the right. The two groups separate symmetrically and nobody change direction.
The turning point trajectories are plotted in white (WFT) and magenta (DPA). In all the simulations we fix the space discretization step in the WFT algorithm as $\varepsilon = 10^{-4}$ and the number of particles $N = 2000$ for the (DPA).}
\label{f:01}
\end{figure}

\begin{figure}[!ht]\centering
\includegraphics[scale=0.357]{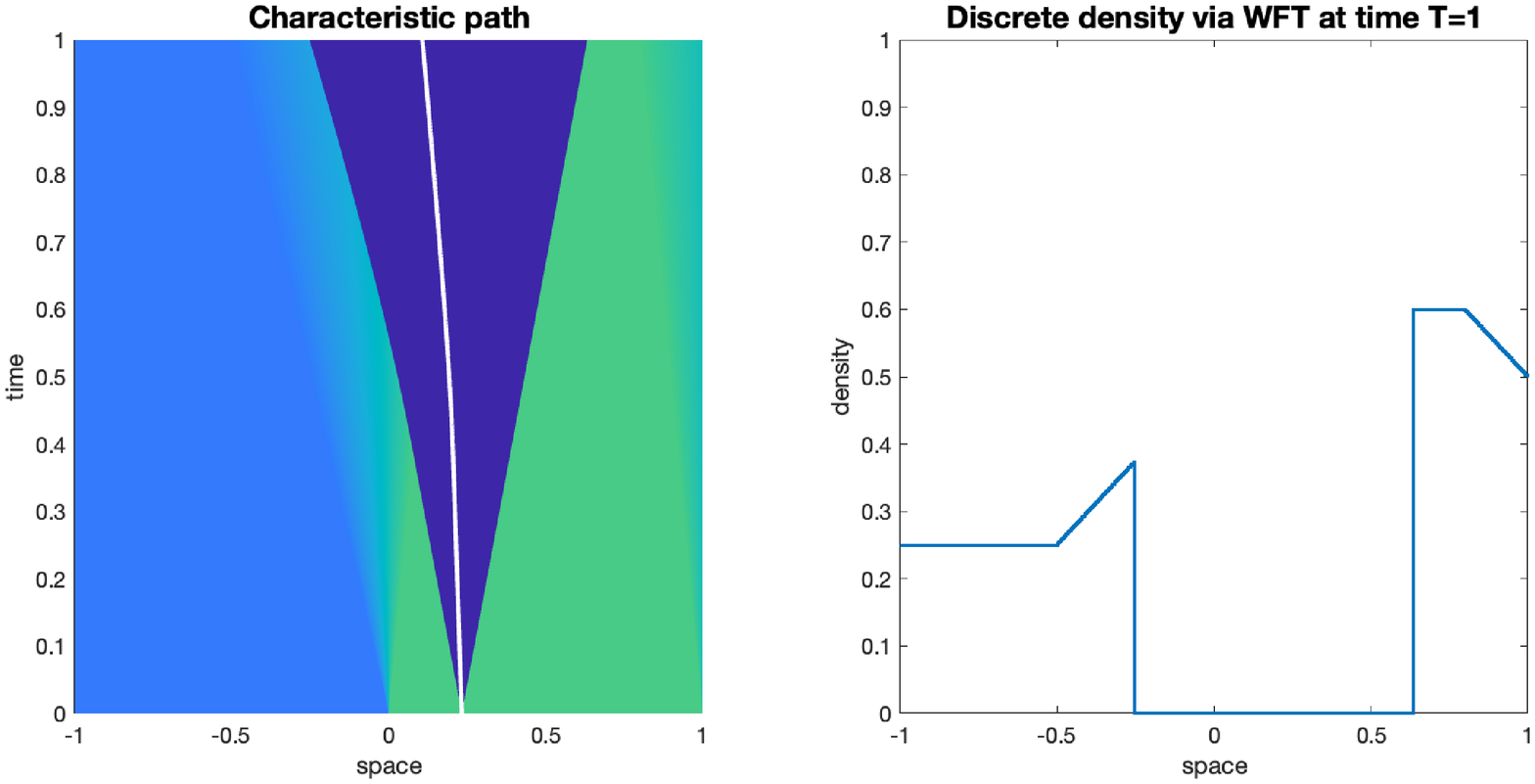}\\
\includegraphics[scale=0.357]{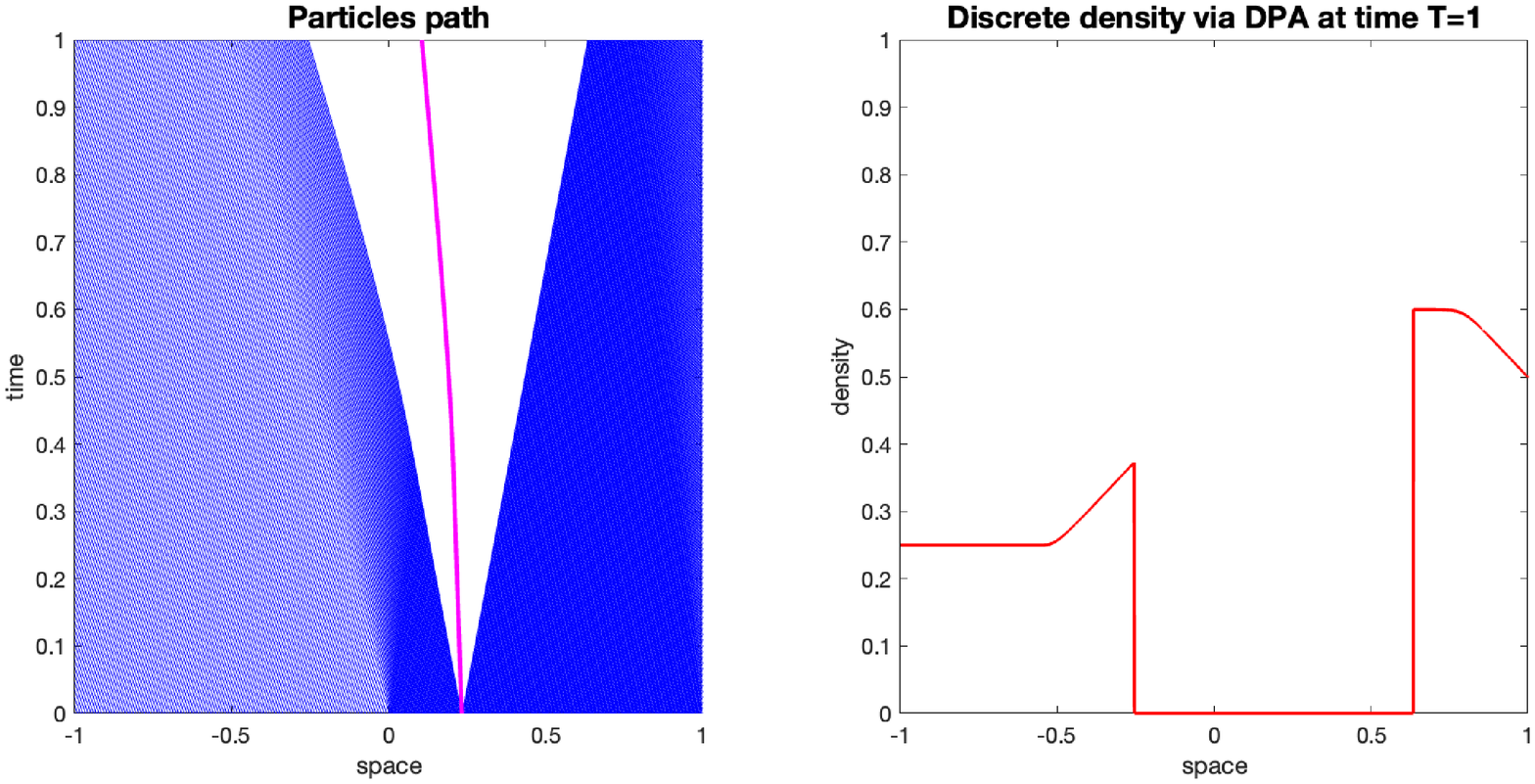}
\caption{
Approximate solutions of Hughes' model with  initial  condition~\eqref{e:infig2} constructed by the WFT scheme (top) and DPA method (bottom). The spatio-temporal evolution is depicted on the left, the corresponding density profile at time $t=1$ on the right. Again, the two groups separate  and nobody change direction.
The turning point trajectories are plotted in white (WFT) and magenta (DPA).}
\label{f:02}
\end{figure}

\begin{figure}[!ht]\centering
\includegraphics[scale=0.357]{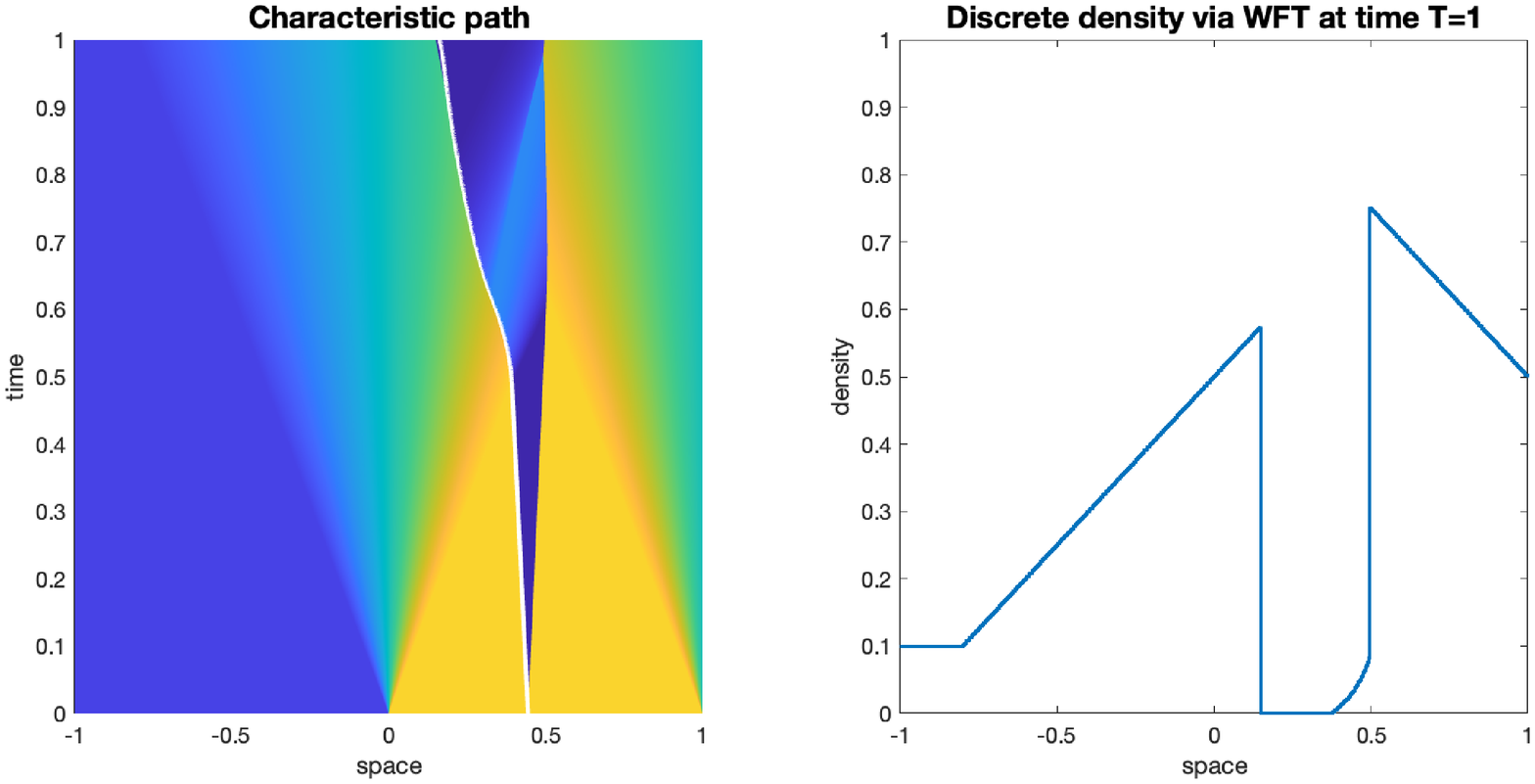}\\
\includegraphics[scale=0.357]{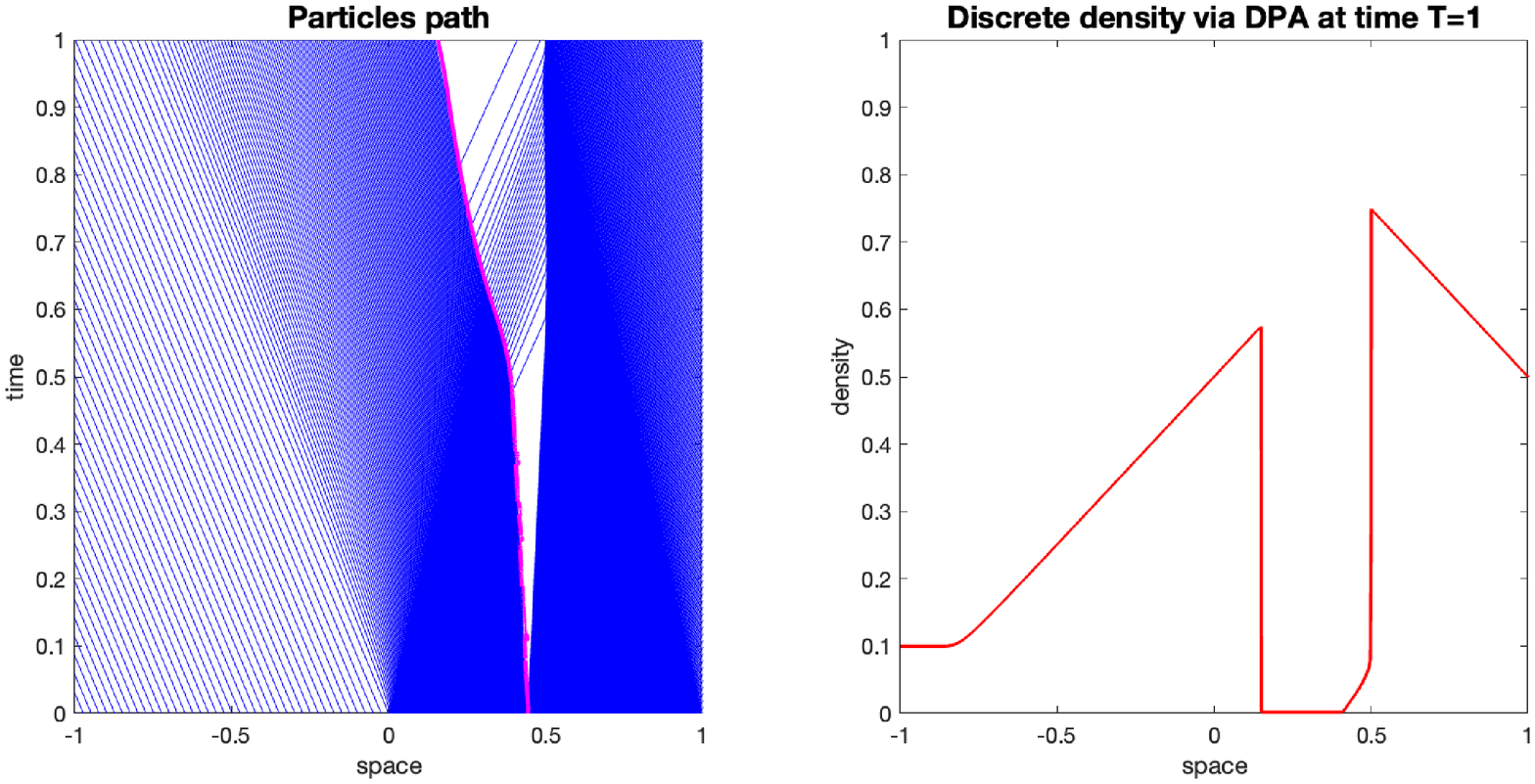}
\caption{
Approximate solutions of Hughes' model with  initial  condition~\eqref{e:infig3} constructed by the WFT scheme (top) and DPA method (bottom). The spatio-temporal evolution is depicted on the left, the corresponding density profile at time $t=1$ on the right. This time, the two groups separate initially, but later some of those initially moving to the left change direction and cross the turning curve.
The turning point trajectories are plotted in white (WFT) and magenta (DPA).}
\label{f:03}
\end{figure}

\begin{figure}[!ht]\centering
\includegraphics[scale=0.357]{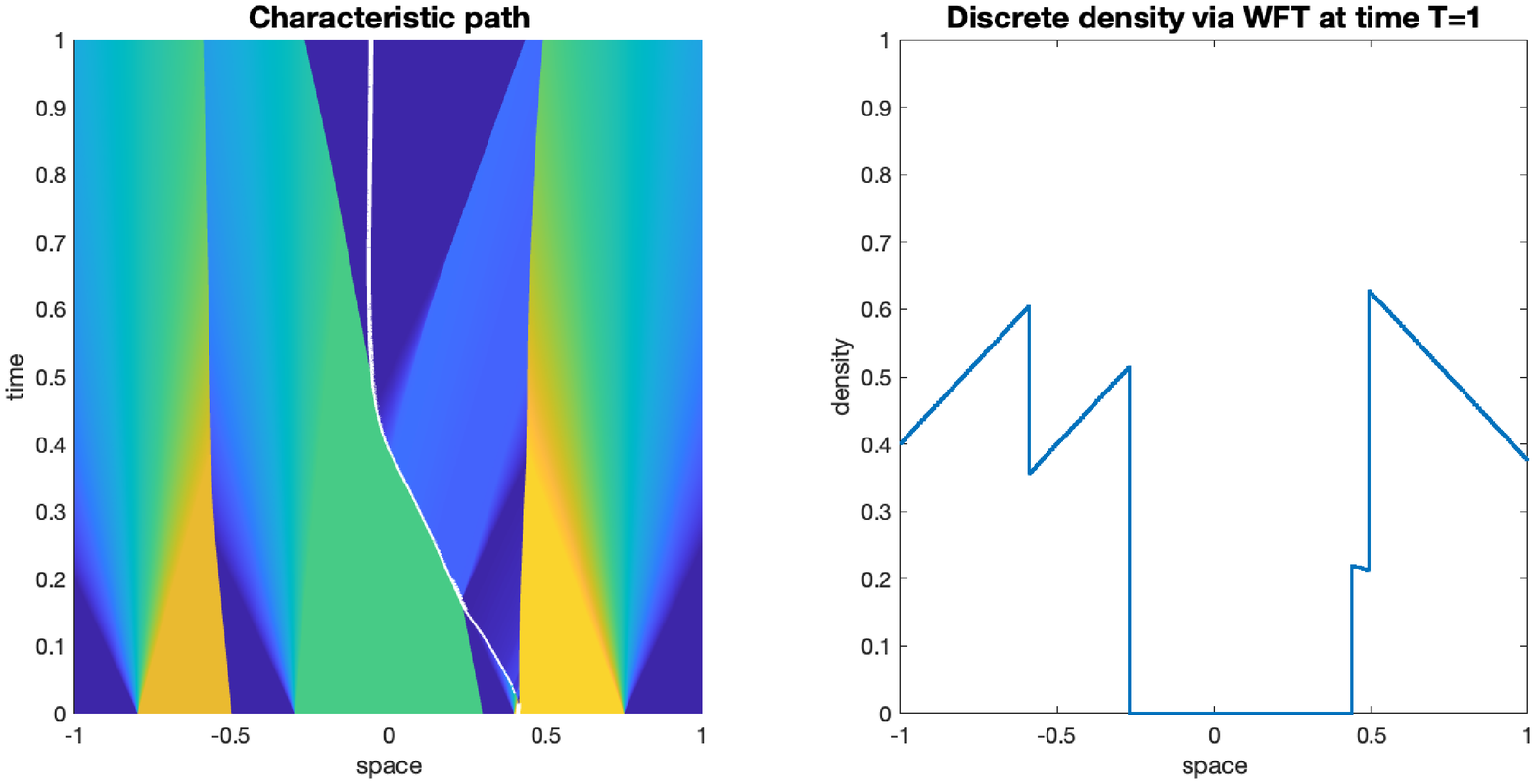}\\
\includegraphics[scale=0.357]{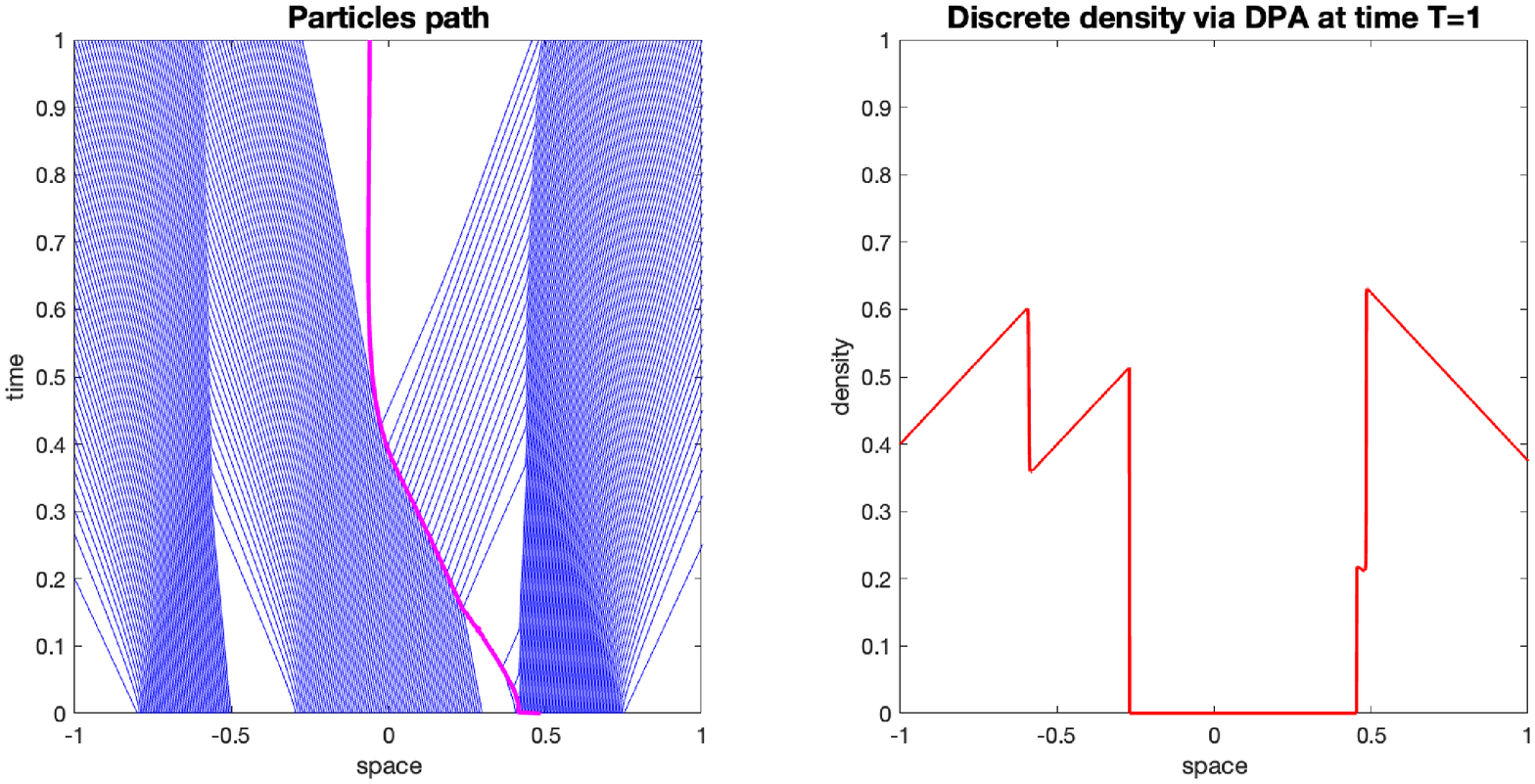}
\caption{
Approximate solutions of Hughes' model with  initial  condition~\eqref{e:infig4} constructed by the WFT scheme (top) and DPA method (bottom). The spatio-temporal evolution is depicted on the left, the corresponding density profile at time $t=1$ on the right. Also here, the two groups separate initially, but later some of those initially moving to the left change direction and cross the turning curve.
The turning point trajectories are plotted in white (WFT) and magenta (DPA).}
\label{f:04}
\end{figure}

\textbf{Wave-Front Tracking (WFT) scheme.} According to the method introduced in \cite{MR3055243}, consider the grid $\mathcal{G}^n$ and the piecewise linear function $f^n$ as introduced in Section \ref{s:WFT}. We approximate the initial datum $\bar{\rho}$ by taking a piecewise constant function
\begin{equation*}
\bar{\rho}^n(x) = \sum_{j} \bar{\rho}_j^n \mathbbm{1}_{(x_{j-1},x_j]}(x),   
\end{equation*}
with density values $\bar{\rho}_j^n\in\mathcal{G}$ and  
jump points $x_j$ such that $\bar{\rho}^n$ approximates $\bar{\rho}$ in the strong $\L1$ topology. Define $\bar{\xi}^n$ according to \eqref{e:barxiWFT}, that reduces to
\begin{equation*}
    \sum_{j\leq 0}c(\bar{\rho}_j^n)(x_j-x_{j-1}) =  \sum_{j> 0}c(\bar{\rho}_j^n)(x_j-x_{j-1}).
\end{equation*}
Then, up to the first waves collision, the approximate solution $\left(\rho^n,\xi^n\right)$ is constructed by solving locally the  Riemann problems
\begin{equation*}
    \begin{dcases}
        \partial_t \rho +\partial_x\left(\operatorname{sign}\left(x-\bar{\xi}^n\right)f^n(\rho)\right) = 0,\\
        \rho(0,x) = \begin{dcases}
            \bar{\rho}_0^n&\mbox{ if }x<\bar{\xi}^n,\\
            \bar{\rho}_1^n&\mbox{ if }x>\bar{\xi}^n,
        \end{dcases}\\
        \dot{\xi}\left(\rho^{+}-\rho^{-}\right) = \Psi[\rho],
    \end{dcases}\quad
        \begin{dcases}
        \partial_t \rho +\partial_x\left(\operatorname{sign}\left(x_j-\bar{\xi}^n\right)f^n(\rho)\right) = 0,\\
        \rho(0,x) = \begin{dcases}
            \bar{\rho}_j^n&\mbox{ if }x<\bar{\xi}^n,\\
            \bar{\rho}_{j+1}^n&\mbox{ if }x>\bar{\xi}^n,
        \end{dcases}\\
       \hfill j\neq 0,
    \end{dcases}
\end{equation*}
where $\Psi[\rho]=\dot{\xi}\left(c(\rho^{+})+c(\rho^{-})\right)$. Note that the solution to the Riemann problem on the left, should be understood by means of the Riemann solver described previously. This procedure generates a new set of values $\rho_{1,j}^n\in\mathcal{G}^n$ and corresponding waves between consecutive densities moving with speeds $\lambda_{1,j}$ determined by the Rankine-Hugoniot condition, that allow the definition of the wave trajectories $x_j(t) = x_j+\lambda_{1,j}t$ and the piecewise constant approximate density
\begin{equation*}
\rho^n(t,x) = \sum_{j} \rho_{1,j}^n \mathbbm{1}_{(x_{j-1}(t),x_j(t)]}(x).
\end{equation*}
Accordingly, a turning curve trajectory can be defined by $\xi^n(t) = \bar{\xi}^n+\dot{\xi}^{n}t$. When two waves collide, a wave hits the domain boundary or the turning curve, new Riemann problems arise, either solved by the classical method or according to the Riemann solver designed in Section~\ref{sec:RP}.

The \textsc{Matlab} code used for the simulations presented in this section can be downloaded at the following URL:\\
\url{http://www-sop.inria.fr/members/Paola.Goatin/wft.html}

\medskip

\textbf{Deterministic Particle Approximation (DPA).}  Concerning the DPA approach, given an initial datum $\bar{\rho}$, we construct the initial particle configuration according to \eqref{e:iniDPA}, which corresponds to a piecewise uniform grid for piecewise constant initial condition, and we solve the particle system \eqref{e:FtL} using the Runge-Kutta \textsc{Matlab} solver ODE23s.  We then reconstruct the density according to 
\begin{equation*}
\rho^n(t,x) = \sum_{i=0}^{N-1} R_{i+1/2}(t) \, \mathbbm{1}_{[\mathrm{x}_i(t),\mathrm{x}_{i+1}(t))}(x),
\end{equation*}
with
\begin{equation*}
R_{i+1/2}(t) = \begin{dcases}
    \frac{m}{\mathrm{x}_{i+1}(t)-\mathrm{x}_i(t)}, & i \in \left\{0,\ldots,N-1\right\}\setminus\left\{I_0\right\},
\\
 0,& i \in \left\{-1,I_0,N\right\}.
\end{dcases}
\end{equation*}
 
An important remark has to be stated about the boundary conditions. As described previously, we do not impose any boundary condition in the particle method. The two leading particles $\mathrm{x}_0$ and $\mathrm{x}_N$ move with maximal velocity towards the opposite directions. 

Particular attention is devoted to the turning point evolution in the particle simulations, obtained by discretizing \eqref{e:xiDPA}. Since no boundary conditions are imposed for the particle method, particles are free to exit the domain following the evolution of the two leaders, whereas only the particles still inside the domain bias the evolution of the turning point. Note that the particles trajectories $\mathrm{x}_i$ constructed through the DPA are in general different from the wave trajectories $x_j(t)$ produced by the WFT algorithm.

\medskip

In all the simulations, we fix the space discretization step in the WFT algorithm as $\varepsilon = 10^{-4}$ and the number of particles $N = 2000$ for the DPA and we plot the reconstructed piecewise constant densities at time $T=1$ and characteristic and particles path for WTF and DPA respectively. In \figurename~\ref{f:01}, and \figurename~\ref{f:02} we plot the numerical solutions corresponding to the initial data
\begin{equation}\label{e:infig1}
    \bar{\rho}(x)=0.6\quad\mbox{ for all }\,x\in\mathsf{D},
\end{equation}
and 
\begin{equation}\label{e:infig2}
    \bar{\rho}(x)=\begin{dcases}
        0.25&\mbox{ if }\,x\in\left[-1,0\right),\\
        0.6&\mbox{ if }\,x\in\left[0,1\right].
    \end{dcases}
\end{equation}
This examples show the classical split of the density in two subgroups moving towards the opposite exits. A more interesting behaviour is shown in 
\figurename~\ref{f:03} and \figurename~\ref{f:04}, were collisions between wave/particles and the turning curve occur, leading to the formation of non-classical shocks. Initial data in  \figurename~\ref{f:03} and \figurename~\ref{f:04} are given respectively by
\begin{equation}\label{e:infig3}
    \bar{\rho}(x)=\begin{dcases}
        0.1&\mbox{ if }\,x\in\left[-1,0\right),\\
        0.9&\mbox{ if }\,x\in\left[0,1\right],
    \end{dcases}
\end{equation}
and
\begin{equation}\label{e:infig4}
    \bar{\rho}(x)=\begin{dcases}
        0.8&\mbox{ if }\,x\in\left[-0.8,-0.5\right),\\
        0.6&\mbox{ if }\,x\in\left[-0.3,0.3\right],\\
         0.9&\mbox{ if }\,x\in\left[0.4,0.75\right],\\
        0&\mbox{ otherwise}\,,
    \end{dcases}
\end{equation}
see also the discussion in  \cite[Section~8]{AndRosSti}.

\section{Modified versions}\label{sec:modified}

We conclude by presenting some generalisations or slight modifications of Hughes' model, which on the one hand try to make it more realistic and on the other hand give insights into the mathematical modelling (especially the missing microscopic interpretation of the original model by Hughes).\\
We discuss three different approaches; we first introduce a \emph{regularised} version of the model proposed in \cite{markowich}, which renders the mathematical theory of the model more accessible through relatively standard techniques; then we describe the  variational approach proposed by Burger et al.~\cite{burger,HerzogPietschmannWinkler} which leads to a dynamic version of Hughes' model; finally, we describe a generalisation of the model proposed by Carrillo et al., see \cite{carrillo}, which is meant to make the model closer to real situations by removing the assumption that each pedestrian has a \emph{global} view of the distribution of the crowd on the whole domain. Further generalisations, including modification of exit behaviour of agents accounting for capacity drop phenomena, were very recently proposed in \cite{AndrGirard-preprint}: the corresponding existence results were sketched in Section~\ref{s:fixed-point}.

\subsection{The regularised Hughes model}
\label{s:regularised}

As mentioned in the introduction, the main difficulty in developing a mathematical theory for the Hughes model resides in the discontinuity of $\nabla \phi$ in the equation \eqref{e:Hughes1D}. 
To bypass this problem, the authors in \cite{markowich} proposed a \emph{regularised version} of the model in one space dimension, in which the eikonal equation \eqref{e:Hughes1D-b} contains extra terms to avoid discontinuities:
\begin{itemize}
\item
an additive constant in the denominator of the right-hand side of the eikonal equation in order to avoid infinite slopes for $\phi$;
 \item 
 a Laplacian term in the (squared) eikonal equation to smoothen the potential $\phi$.
\end{itemize}
The resulting model considered in \cite{markowich} is
\begin{subequations}\label{eq:firstreg}
\begin{empheq}[left=\empheqlbrace]{align}
\label{eq:firstreg_rho}
    &\rho_t- \bigl(\rho v(\rho)^2 \phi_x\bigr)_x =0, \\
    &-\delta_1\phi_{xx} + |\phi_x|^2 = \frac{1}{\bigl(v(\rho)+\delta_2\bigr)^2},
\end{empheq}
\end{subequations}
where $\delta_1,\delta_2>0$ are two parameters of the model.
Therefore, the result in \cite{markowich} assumes $g\equiv1$. 
Moreover, it is assumed that $v(\rho)=(1-\rho)_+$ for simplicity. 
The model is posed on a bounded interval $x\in \mathsf{C} \doteq (-1,1)$ with homogeneous Dirichlet boundary conditions
\begin{align*}
&\rho(\pm 1^{\mp},t)=0\,,&&\phi(\pm1^{\mp},t)=0\,.
\end{align*}

A suitable notion of entropy solution can be formulated for \eqref{eq:firstreg}, where the modified eikonal equation is solved with respect to $\phi$, $\phi=\phi[\rho](t,x)$. Clearly, such a dependence is \emph{non-local} in space.  This notion of entropy solutions is a natural generalisation of cases previously considered in the literature of scalar conservation laws with space-time dependent fluxes. Roughly speaking, it is assumed:
\begin{itemize}
\item that $\rho$ is continuous in time with values in $\BV([-1,1])$;
\item that $\phi$ is continuous in time with values in $\W{2,\infty}$;
\item that $\rho$ and $\phi$ satisfy the entropy inequality
\end{itemize}
\begin{align}
&\int_0^{\infty}\int_{-1}^1 \Bigl( |\rho-k| \, \varphi_t + \mathrm{sign}(\rho-k)m(k) \bigl( \varphi \, \phi_{xx} - \bigl(m(\rho)-m(k)\bigr) \, \phi_x \, \varphi_x \bigr) \Bigr) \,\d x\,\d t
\nonumber\\ 
&-\mathrm{sign}(k)\int_0^T \bigl(m(\mathrm{tr}(\rho))-m(k)\bigr) \, \phi_x \, \varphi\Big\vert_{x=\pm 1}\,\d t
+\int_{-1}^1 \bar\rho(x) \, \varphi(0,x)\,\d x
\ge 0,\label{eq:entropysolution}
\end{align}
\begin{itemize}
\item[]
where $m(\rho) \doteq \rho v(\rho)^2$ 
and $\varphi$ is an arbitrary $\Cc\infty$ test function;
\item that $\rho$ and $\phi$ satisfy the regularised eikonal equation in \eqref{eq:firstreg} almost everywhere.
\end{itemize}
The above entropy inequality \eqref{eq:entropysolution} incorporates the concept of entropy solutions at the boundary according to the classical approach in \cite{BardosLerouxNedelec}. We refer to \cite{markowich} for details. 

The existence of entropy solutions according to the above definition is carried out in \cite{markowich} by a standard vanishing viscosity approach for the continuity equation. A major issue to achieve the needed compactness is the regularity of the potential $\phi$. To perform this task, the authors apply a Hopf-Cole transformation
\[\psi(x,t) \doteq e^{-\frac{\phi(x,t)}{\delta_1}},\]
which implies the following boundary value problem for $\psi$
\[\begin{cases}
\displaystyle{\delta_1^2 \psi_{xx} =\psi F_{\delta_2}(\rho)} \,,\\
    \psi(\pm 1)= 1\,, &
\end{cases}\]
where
\[F_{\delta_2}(\rho) \doteq \frac{1}{(v(\rho)+\delta_2)^2}\,.\]
Then, a simple multiplication by $\psi$ and integration by parts on $[-1,1]$ imply $\W{2,\infty}$ regularity for $\psi$ and, consequently, for $\phi$. We refer to \cite[Lemmas~2.2. and~2.3]{markowich}. 
These results allow to obtain uniform bounds in $\L\infty\cap \BV$ for $\rho$ with respect to the artificial viscosity parameter and to obtain convergence up to a subsequence to an entropy solution. The uniqueness is obtained by standard doubling of the variables. We refer to \cite{markowich} for further details.

\subsection{A dynamic version of Hughes model via optimal control}

The dynamic formulation is based on the assumption that pedestrians wish to exit a domain at fastest. This corresponds to a classical or stochastic optimal control problem on the microscopic level, and a PDE constrained optimisation problem on the macroscopic level. We will see that we can relate the optimality conditions of this transient optimal control problem to Hughes' model in suitable scaling limits. 

Consider a pedestrian (of unit mass) trying to leave the domain $D \subset \mathbb{R}^2$ as fast as possible. Let $\mathbf{X} = \mathbf{X}(t)$ denote its position at time $t>0$, $\mathbf{V} = \mathbf{V}(t)$ its velocity and $\mathbf{X}_0$ its starting position. Define the exit time as 
$$ T_{\rm exit}(\mathbf{X}) = \sup\{ t > 0 : \mathbf{X}(t) \in D\}. $$
We assume that pedestrians are perfectly rational and wish to minimise a weighted sum of the exit time $T_{\rm exit}$ and the kinetic energy, i.e.
\begin{equation}\label{eq:energy1}
	\frac{1}2 \int_0^{T_{\rm exit}} |\mathbf{V}(t)|^2\,\d t + \frac{\alpha}2 T_{\rm exit}(\mathbf{X}) \rightarrow \min_{\mathbf{X},\mathbf{V}} ,
\end{equation}
subject to $\dot{\mathbf{X}}(t) = \mathbf{V}(t)$, $\mathbf{X}(0) = \mathbf{X}_0$ and given a weighting parameter $\alpha > 0$.  Next we introduce the Dirac measure $\mu = \delta_{\mathbf{X}(t)},$ and choose the final time $T$ sufficiently large. Then we can rewrite \eqref{eq:energy1} in terms of the measure $\mu$:
\begin{equation}\label{eq:energy2}
 I_T(\mu,v) = \frac{1}2 \int_0^T \int_D|v(t,x)|^2 \d\mu\,\d t + \frac{\alpha}2 \int_0^T \int_D \d\mu\,\d t ,
\end{equation}
subject to the constraint $\mu_t + \nabla \cdot (\mu v) = 0$,
with initial condition $\mu|_{t=0} = \delta_{\mathbf{X}_0}$. By reformulating \eqref{eq:energy1} in terms of the measure $\mu$ we obtain the continuum version \eqref{eq:energy2}. \\
If we loosen the rationality assumption and allow for uncertainty in the pedestrian's path the ODE for $\mathbf{X}$ is replaced by 
a stochastic differential equation
\begin{subequations}\label{e:stochopt}
\begin{equation}
	\d\mathbf{X}(t) = \mathbf{V}(t)\,\d t + \sigma\,\d W(t), \label{Langevin}
\end{equation}
where $W$ is a Wiener process and $\sigma$ the diffusivity. Due to the stochasticity we consider the expected value of \eqref{eq:energy1}, that is
\begin{equation}
\mathbb{E}_{\mathbf{X}_0} \left[	\frac{1}2 \int_0^{T_{\rm exit}} |\mathbf{V}(t)|^2\,\d t + \frac{\alpha}2 T_{\rm exit}(\mathbf{X})
\right] \rightarrow \min_{\mathbf{V}}
\end{equation}
\end{subequations}
with the random variable $\mathbf{X}$ determined by \eqref{Langevin} with initial value $\mathbf{X}_0$. \\
Rewriting \eqref{e:stochopt} in terms of the distribution $\mu$ and assuming that $\mu$ has a density $\rho$, that is $\d\mu = \rho \,\d x$, gives the respective macroscopic formulation:
\begin{equation*}
 I_T(\rho,v) = 	\frac{1}2 \int_0^T \int_D \rho(t,x)~ |v(t,x)|^2 \,\d x\,\d t + \frac{1}2 \int_0^T \int_D \rho(t,x)\,\d x\,\d t ,
\end{equation*}
subject to $\rho_t + \nabla \cdot (\rho v) = \frac{\sigma^2}2 \Delta \rho,$  with $\rho(x,0) = \rho_0(x).$
The formal optimality conditions of this constrained optimisation problem are
\begin{equation}
\label{e:optimality}
\begin{cases}
	\displaystyle\rho_t + \nabla \cdot ( \rho \nabla \phi ) - \frac{1}{2} \sigma^2 \Delta \rho = 0, \\
	\displaystyle\phi_t + \frac{1}{2} \|\nabla \phi\|^2  + \frac{1}{2} \sigma^2 \Delta \phi = \frac{\alpha}2,
\end{cases}
\end{equation}
where $\phi$ corresponds to the dual or adjoint variable. Note that system \eqref{e:optimality} is supplemented with an initial condition for $\rho(x,0) = \rho_0(x)$ and a terminal condition for $\phi(x,T) = 0$. The adjoint variable $\phi$ satisfies a transient viscous eikonal equation, which has to be solved backward in time. The connection to the Hughes model is quite apparent for a time interval $[0,S]$ with $S \ll T$ and $\sigma = 0$. Noticing that the Hamilton-Jacobi equation is solved backwards in time and that the backward time is large for $t \leq S$, we see that  the solution $\phi$ is mainly determined by the large-time asymptotics solving, for some $c \in \mathbb{R}^+$,
\[\|\nabla \tilde \phi\|^2 = c.\]

\noindent Motivated by the above interpretation of the Hughes problem Burger et al.\ investigated the following generalisation on the macroscopic level
\begin{subequations}\label{e:optcon}
\begin{equation} \label{e:opt}
  I_T(\rho,v) \doteq	\frac{1}2 \int_0^T \int_D F(\rho) \, |v(t,x)|^2 \,\d x\,\d t + \frac{1}2 \int_0^T \int_D E(\rho)\,\d x\,\d t ,
\end{equation}
subject to
\begin{equation}\label{e:constraint}
	\rho_t + \nabla \cdot \bigl(G(\rho) v\bigr) = \frac{\sigma^2}2 \, \Delta \rho,
\end{equation}
\end{subequations}
and a given initial value $\rho(x,0) = \rho_0(x)$. Here the functions $G$, $F$ and $E$ account for nonlinear effects in high density regimes. In particular
\begin{itemize}
\item The function $G = G(\rho)$ corresponds to a nonlinear mobility. In the setting of pedestrian dynamics it is assumed to be a positive non-negative function of the density. For example, $G$ is often assumed to be non-increasing and approaching zero when approaching the maximum capacity.
\item The function $F=F(\rho)$ accounts for the modulation of transport costs by the density. For example the function $F$ might tend to infinity as $\rho$ approaches $\rho_{\max}$.
\item The nonlinear function $E=E(\rho)$ in the exit time functional can for example relate to increased cost of moving in high density regions.
\end{itemize}
\medskip

Burger et al.\ discussed the relation of the original Hughes model to the solution of the optimality system defined by \eqref{e:optcon}. They showed that for 
\begin{align*}
&\sigma = 0,&
&H(\rho) = \frac{G^2(\rho)}{F(\rho)} = \rho f(\rho)&
&\text{ and }&
&E(\rho) = \alpha \rho,
\end{align*}
with $f(\rho) = \rho_{\max} - \rho$, the optimality system of \eqref{e:optcon} is given by
\begin{subequations}
\begin{align}
	\rho_t + \nabla \cdot \bigl( \rho~f(\rho)^2 \nabla \phi \bigr) &= 0 \\
	\phi_t + \frac{f(\rho)}2 \bigl(f(\rho) + 2 \rho f'(\rho)\bigr) \|\nabla \phi\|^2  &= \frac{\alpha}2.\label{e:modhughphi}
\end{align}
\end{subequations}
Arguing again that the Hamilton-Jacobi equation equilibrates much faster for large times $T$ and therefore neglecting $\phi_t$ in \eqref{e:modhughphi} gives a problem which is quite close to the original Hughes model for $f(\rho) = \rho_{\max} - \rho$, but still, with a different prefactor in the eikonal equation (due to the $2\rho f'(\rho)$ term).

 Therefore Burger et al.\ provided another formal argument which links the original Hughes model to the dynamic formulation. They proposed a modified mean field approach by extrapolating the current density into the future, 
i.e.~$\rho = \rho(x,t)$ is assumed to be the density of the system for all times $s > t $. 

\noindent In doing so they consider $N$ particles with position $\mathbf{X}_k = \mathbf{X}_k(t)$ and define the empirical density 
\begin{align*}
\rho^N(t) =  \frac{1}{N} \sum_{k=1}^N \delta\bigl(x-\mathbf{X}_k(t)\bigr).
\end{align*}
Furthermore they introduce a smoothed approximation $\rho^N_\zeta$ of the empirical density (which is necessary to define the cost functional later)
\[
\rho^N_\zeta(t) = (\rho^N * \zeta)(t,x) = \frac{1}{N} \sum_{k=1}^N  \zeta\bigl(x-\mathbf{X}_k(t)\bigr),
\]
for a sufficiently smooth positive kernel $\zeta$. \\ 
Based on the considerations above they assume that the optimal velocity of an agent at position 
$\mathbf{X} = \mathbf{X}(t)$ is determined by minimising
\begin{equation}\label{e:minC}
 \frac{1}2\int_t^{t+T} \frac{\lvert \mathbf{V}(s)\rvert ^2}{G\bigl(\rho_\zeta^N(\xi(s;t),t\bigr)} \,\d s + \frac{1}2 T_{\rm exit}(\mathbf{X},\mathbf{V}) \rightarrow  \min_{(\mathbf{X},\mathbf{V})}
\end{equation}
subject to the constraint that $\frac{\d \xi}{\d s}=\mathbf{V}(s)$ and $\xi(0)=\mathbf{X}(t)$. Hence an individual tries to find its optimal trajectory based on the current density $\rho$ and extrapolating it into the future. We can rewrite the above problem in terms of probability measures (replacing the smoothed empirical density $\rho^N_\zeta$ by its mean field limit $\rho$) and obtain
\[
	\frac{1}2 \int_t^{T+t}  \int_\Omega \left(\frac{w^2(x,s)}{G\bigl(\rho(\xi(s;t),t)\bigr)} + 1 \right)\,\d \mu \,\d s \rightarrow \min_{(\mu,w)}
\]
for the velocity field $w$ and the probability measure $\mu$  satisfying $\mu_s + \nabla \cdot(\mu w) = 0.$ with initial condition $ \mu(t=0)=\delta_X$.
Formal calculation of the optimality conditions for $G(\rho) = f(\rho) = \rho_{\max} - \rho$, and the argument that the adjoint variable equilibrates faster for sufficiently large $T$ then gives the original Hughes model. Hence we can interpret Hughes model as a microscopic optimal control problem, in which agents determine the optimal trajectory using the current pedestrian density and interpolating it into the future. 

\subsection{Optimal control via local attraction}

Related is the optimal control problem discussed in \cite{HerzogPietschmannWinkler} which is based on a regularised version similar to the one illustrated in Section \ref{s:regularised}, yet with an additional Laplacian term added in \eqref{eq:firstreg_rho}. The idea is to control the trajectories of a fixed, finite number $M$ of agents that are able to influence the crowd in the vicinity of their location. A typical example, thinking about a tourist guide or security personal, would be an local, attractive force. This is included into the model by an additional convection term in \eqref{eq:firstreg_rho} which is the gradient of an interaction kernel centred at the agents location (denoted by $x_i(t)$, $i\in\{1,\ldots, M\}$). Thus \eqref{eq:firstreg_rho} becomes
\begin{align*}
\rho_t- \nabla \cdot \Biggl(\rho v(\rho)^2 \biggl( \nabla \phi + \sum_{i=1}^M \nabla K\bigl(x-x_i(t)\bigr)\biggr)\Biggr) =\delta_3 \Delta \rho,
\end{align*}
with an attractive interaction kernel $K$, typically radially symmetric and with compact support, and $\delta_3 > 0$.
The motion of the agents themselves is then governed by an ordinary differential equation of the form 
\begin{align}\label{eq:firstreg_rho_reg}
&\dot{x}_i(t) = v\bigl(\rho(x_i(t),t)\bigr) u_i(t),&
&x_i(0) = x_i^0,&
&i\in\{1,\ldots, M\}.
\end{align}
The vector fields $u_i$ are the actual controls that determine the agent's direction while the first term on the right hand side makes sure that the agents are slowed in high density areas as is the remaining crowd. However, as this requires a point evaluation of $\rho$ at $x_i(t)$, sufficient regularity needs to be shown which is the reason for the additional diffusive term in \eqref{eq:firstreg_rho_reg}. Possible objective functionals are the total mass at some final time (evacuation scenario) or the area of parts of the domain in which a given critical density is exceeded (panic avoidance). Control of crowds via few agents has also been considered in \cite{AlbiBonginiCristianiKalise:2016:1, BurgerPinnauRothTotzeckTse:2016:1}, yet for different models and different applications. The main result of \cite{HerzogPietschmannWinkler} is the existence of the regularised model with sufficient regularity as well as differentiability properties of the control-to-state map, existence of a globally optimal control, and the formulation optimality conditions. In a subsequent work, \cite{PietschmannStoetznerWinkler}, the authors introduce a numerical discretisation based on a finite volume scheme that is shown to preserve the box constraints of $\rho$ and provide a variety of numerical examples of the optimal control problem.

\subsection{A localised version of the model}

In the original Hughes model the function $\phi$ is calculated assuming that the global distribution of pedestrians is known at any time $t>0.$ The assumption of global knowledge of the pedestrian density is highly questionable in practical situations, which is why Carrillo et al.~\cite{carrillo} proposed a local version to account for limited vision and restricted perception. Their starting point was motivated by the microscopic interpretation of Hughes model \eqref{e:minC}. However, they assumed that individuals can only estimate the pedestrian density within their vision cone. Since \eqref{e:minC} is a first order model, the implementation of a vision cone is not as straightforward than in a second order model. \\
Therefore Carrillo et al.~introduced an auxiliary variable and a parametrised potential $\phi(x,y) \colon \mathbb{R}^4\rightarrow\mathbb{R}$ such that $y\mapsto \phi(x_0,y)$ denotes the cost potential calculated by pedestrians located at $x_0\in\Omega$. For every point $x$ we assume that the domain $D$ is decomposed into a visible part $D_V$ and an
invisible one  $D_I = D \backslash D_V$.
Then the limited perception can be implemented as follows: 
in visible areas the optimal trajectory is calculated using the pedestrian
density, while in invisible areas the density is set to a constant 
value $\rho_H \in \R^+_0$. We assume that $\rho_H$ is the same for all pedestrians. 
For example, if $\rho_H=0$ then pedestrians assume that not visible areas are
empty, while pedestrians will avoid
these areas if $\rho_H\approx \rho_{\max}$. 
The respective eikonal equation is then
\[
\lVert \nabla_y \phi(x,y )\rVert = \displaystyle{\begin{cases}
 1/v\bigl(\rho(y,t)\bigr), & y\in D_V, \\
 1/v(\rho_H), & y\in D_I,
 \end{cases}}
\]
which gives the potential $\phi$ as function of two space
variables. \\
Carrillo et al.\ calculate this potential for every single exit (since the visible and invisible areas change for each one). The final walking direction in each point is then computed by comparing the potentials for all exits and adjusting it according to the predominant direction in the close surrounding. This averaging is necessary to avoid strong fluctuations in the walking direction.  We omit the details of the full model as it exceeds the scope of this review. Computational experiments show that this generalisation yields more realistic results especially in the case of obstacles and more complicated geometries.

\section{Conclusions and future challenges}

As was to be expected given the mathematical structure of the model, the various (successful and unsuccessful) attempts to prove existence of solutions to Hughes' model involve the community of researchers working on hyperbolic conservation laws. As for the one-dimensional case, the shock structure of the model is by now quite well understood. It is somehow surprising, though, that the only existence results for large data do not use the WFT algorithm. Therefore, a first open question is the convergence of the WFT scheme, at least in the case of linear cost, that is the one covered in the available existence theorems. The next step towards a satisfactory one-dimensional theory is to prove existence of entropy solutions for more general cost functionals. Having three approaches which lead to significant results so far, namely the WFT scheme, the DPA scheme, and the fixed point strategy, makes us quite optimistic that this result is within reach.

Having developed significant results in the one-dimensional case was a necessary intermediate step to the solution of  this model and to better understand its mathematical features. However, the journey towards a satisfactory mathematical theory for Hughes' model cannot be considered as completed unless some results are obtained in two space dimensions. There are many possible directions to take in this sense:
\begin{itemize}
    \item Extend the available results on the regularised model to the two-dimensional case. This seems quite reasonable. The strategy adopted so far used the very specific features of the one-dimensional case, but we believe something can be done also in $2d$.
    \item Set up a reasonable deterministic particle scheme in the two dimensional space, for example by using Voronoi tessellation to reconstruct the density. A major issue in this case is the definition of the direction of the discretized flux.
    \item Try to investigate better the structure of coupling with the eikonal equation suitably involving viscosity solutions. In this sense, the interaction with researchers from the field of viscosity solutions should be definitely improved.
\end{itemize}

Models with more general and possibly more realistic boundary conditions need to be further investigated as well, both numerically and analytically, especially in two space dimensions.
More broadly, the interplay with control theory (partly mentioned here) is an almost unexplored direction of research, which we believe would have a relevant impact on the applications and which would certainly benefit from a sound, well-established mathematical theory for the IBV problems.

\subsubsection*{Acknowledgements}
Amadori, Di Francesco, Fagioli, Rosini and Stivaletta are members of GNAMPA-INdAM (\textit{Gruppo Nazionale per l’Analisi Matematica, la Probabilit\`a e le loro Applicazioni - Istituto Nazionale di Alta Matematica}), Russo is a member of GNCS-INDAM (\textit{Gruppo Nazionale per il Calcolo Scientifico}).
Andreianov and Girard would like to thank l’Agence Nationale de la Recherche (ANR) to support this research with funds coming from project ANR-22-CE40-0010 (ANR CoSS). 
Russo would like to thank the Italian Ministry of Instruction, University and Research (MIUR) to support this research with funds coming from PRIN Project 2017 (No. 2017KKJP4X entitled “Innovative numerical methods for evolutionary partial differential equations and applications”).


\begin{thebibliography}{10}
	
	\bibitem{AlbiBonginiCristianiKalise:2016:1}
	G.~Albi, M.~Bongini, E.~Cristiani, and D.~Kalise.
	\newblock {I}nvisible control of self-organizing agents leaving unknown
	environments.
	\newblock {\em SIAM J. Appl. Math.}, 76(4):1683--1710, 2016.
	
	\bibitem{AmadoriDiFrancesco}
	D.~Amadori and M.~Di~Francesco.
	\newblock The one-dimensional {H}ughes model for pedestrian flow:
	{R}iemann-type solutions.
	\newblock {\em Acta Math. Sci. Ser. B (Engl. Ed.)}, 32(1):259--280, 2012.
	
	\bibitem{AmadoriGoatinRosini}
	D.~Amadori, P.~Goatin, and M.~D. Rosini.
	\newblock Existence results for {H}ughes' model for pedestrian flows.
	\newblock {\em J. Math. Anal. Appl.}, 420(1):387--406, 2014.
	
	\bibitem{ADRR2016}
	B.~Andreianov, C.~Donadello, U.~Razafison, and M.~D. Rosini.
	\newblock Qualitative behaviour and numerical approximation of solutions to
	conservation laws with non-local point constraints on the flux and modeling
	of crowd dynamics at the bottlenecks.
	\newblock {\em ESAIM Math. Model. Numer. Anal.}, 50(5):1269--1287, 2016.
	
	\bibitem{ADR2014}
	B.~Andreianov, C.~Donadello, and M.~D. Rosini.
	\newblock Crowd dynamics and conservation laws with nonlocal constraints and
	capacity drop.
	\newblock {\em Math. Models Meth. Appl. Sci.}, 24:2685--2722, 2014.
	
	\bibitem{AndrGirard-preprint}
	B.~Andreianov and T.~Girard.
	\newblock Existence of solutions to a class of one-dimensional models for
	pedestrian evacuations.
	\newblock HAL preprint, https://hal.science/hal-03937464, 2023.
	
	\bibitem{AndRosSti}
	B.~Andreianov, M.~Rosini, and G.~Stivaletta.
	\newblock On existence, stability and many-particle approximation of solutions
	of 1{D} {H}ughes model with linear costs, 2021.
	
	\bibitem{AndrSylla-FVCA}
	B.~Andreianov and A.~Sylla.
	\newblock A macroscopic model to reproduce self-organization at bottlenecks.
	\newblock In {\em Finite volumes for complex applications {IX}---methods,
		theoretical aspects, examples---{FVCA} 9, {B}ergen, {N}orway, {J}une 2020},
	volume 323 of {\em Springer Proc. Math. Stat.}, pages 243--254. Springer,
	Cham, [2020] \copyright 2020.
	
	\bibitem{MR4059365}
	B.~Aylaj, N.~Bellomo, L.~Gibelli, and A.~Reali.
	\newblock A unified multiscale vision of behavioral crowds.
	\newblock {\em Math. Models Methods Appl. Sci.}, 30(1):1--22, 2020.
	
	\bibitem{BardosLerouxNedelec}
	C.~Bardos, A.~Y. le~Roux, and J.-C. N{\'e}d{\'e}lec.
	\newblock First order quasilinear equations with boundary conditions.
	\newblock {\em Comm. Partial Differential Equations}, 4(9):1017--1034, 1979.
	
	\bibitem{MR3149318}
	R.~Borsche, A.~Klar, S.~K\"{u}hn, and A.~Meurer.
	\newblock Coupling traffic flow networks to pedestrian motion.
	\newblock {\em Math. Models Methods Appl. Sci.}, 24(2):359--380, 2014.
	
	\bibitem{MR3177723}
	R.~Borsche and A.~Meurer.
	\newblock Interaction of road networks and pedestrian motion at crosswalks.
	\newblock {\em Discrete Contin. Dyn. Syst. Ser. S}, 7(3):363--377, 2014.
	
	\bibitem{Bressanbook}
	A.~Bressan.
	\newblock {\em Hyperbolic systems of conservation laws}, volume~20 of {\em
		Oxford Lecture Series in Mathematics and its Applications}.
	\newblock Oxford University Press, Oxford, 2000.
	\newblock The one-dimensional Cauchy problem.
	
	\bibitem{burger}
	M.~Burger, M.~Di~Francesco, P.~A. Markowich, and M.-T. Wolfram.
	\newblock Mean field games with nonlinear mobilities in pedestrian dynamics.
	\newblock {\em Discrete Contin. Dyn. Syst. Ser. B}, 19(5):1311--1333, 2014.
	
	\bibitem{BurgerPinnauRothTotzeckTse:2016:1}
	M.~Burger, R.~Pinnau, A.~Roth, C.~Totzeck, and O.~Tse.
	\newblock {C}ontrolling a self-organizing system of individuals guided by a few
	external agents -- particle description and mean-field limit.
	\newblock {arXiv}: \href{https://arxiv.org/abs/1610.01325}{1610.01325}, 2016.
	
	\bibitem{MR3619091}
	F.~Camilli, A.~Festa, and S.~Tozza.
	\newblock A discrete {H}ughes model for pedestrian flow on graphs.
	\newblock {\em Netw. Heterog. Media}, 12(1):93--112, 2017.
	
	\bibitem{MR3698447}
	E.~Carlini, A.~Festa, F.~J. Silva, and M.-T. Wolfram.
	\newblock A semi-{L}agrangian scheme for a modified version of the {H}ughes'
	model for pedestrian flow.
	\newblock {\em Dyn. Games Appl.}, 7(4):683--705, 2017.
	
	\bibitem{carrillo}
	J.~A. Carrillo, S.~Martin, and M.-T. Wolfram.
	\newblock An improved version of the {H}ughes model for pedestrian flow.
	\newblock {\em Math. Models Methods Appl. Sci.}, 26(4):671--697, 2016.
	
	\bibitem{MR3823842}
	R.~M. Colombo, M.~Gokieli, and M.~D. Rosini.
	\newblock Modeling crowd dynamics through hyperbolic-elliptic equations.
	\newblock In {\em Non-linear partial differential equations, mathematical
		physics, and stochastic analysis}, EMS Ser. Congr. Rep., pages 111--128. Eur.
	Math. Soc., Z\"{u}rich, 2018.
	
	\bibitem{ColomboRossi}
	R.~M. Colombo and E.~Rossi.
	\newblock On the micro-macro limit in traffic flow.
	\newblock {\em Rend. Semin. Mat. Univ. Padova}, 131:217--235, 2014.
	
	\bibitem{DafermosWFT}
	C.~M. Dafermos.
	\newblock Polygonal approximations of solutions of the initial value problem
	for a conservation law.
	\newblock {\em J. Math. Anal. Appl.}, 38:33--41, 1972.
	
	\bibitem{Dafermosbook}
	C.~M. Dafermos.
	\newblock {\em Hyperbolic conservation laws in continuum physics}, volume 325
	of {\em Grundlehren der mathematischen Wissenschaften [Fundamental Principles
		of Mathematical Sciences]}.
	\newblock Springer-Verlag, Berlin, fourth edition, 2016.
	
	\bibitem{DiFrancescoFagioliRosini-BUMI}
	M.~Di~Francesco, S.~Fagioli, and M.~Rosini.
	\newblock Deterministic particle approximation of scalar conservation laws.
	\newblock {\em Boll. Unione Mat. Ital.}, 10(3):487--501, 2017.
	
	\bibitem{MR3644595}
	M.~Di~Francesco, S.~Fagioli, M.~Rosini, and G.~Russo.
	\newblock Follow-the-leader approximations of macroscopic models for vehicular
	and pedestrian flows.
	\newblock In {\em Active particles. {V}ol. 1. {A}dvances in theory, models, and
		applications}, Model. Simul. Sci. Eng. Technol., pages 333--378.
	Birkh\"{a}user/Springer, Cham, 2017.
	
	\bibitem{DiFrancescoFagioliRosiniRussoKRM}
	M.~Di~Francesco, S.~Fagioli, M.~D. Rosini, and G.~Russo.
	\newblock Deterministic particle approximation of the {H}ughes model in one
	space dimension.
	\newblock {\em Kinet. Relat. Models}, 10(1):215--237, 2017.
	
	\bibitem{DiFrancescoFagioliRosiniRusso}
	M.~Di~Francesco, S.~Fagioli, M.~D. Rosini, and G.~Russo.
	\newblock A deterministic particle approximation for non-linear conservation
	laws.
	\newblock In {\em Theory, numerics and applications of hyperbolic problems.
		{I}}, volume 236 of {\em Springer Proc. Math. Stat.}, pages 487--499.
	Springer, Cham, 2018.
	
	\bibitem{markowich}
	M.~Di~Francesco, P.~A. Markowich, J.-F. Pietschmann, and M.-T. Wolfram.
	\newblock On the {H}ughes' model for pedestrian flow: the one-dimensional case.
	\newblock {\em J. Differential Equations}, 250(3):1334--1362, 2011.
	
	\bibitem{DiFrancescoRosini}
	M.~Di~Francesco and M.~Rosini.
	\newblock Rigorous derivation of nonlinear scalar conservation laws from
	follow-the-leader type models via many particle limit.
	\newblock {\em Arch. Ration. Mech. Anal.}, 217(3):831--871, 2015.
	
	\bibitem{MR4026959}
	M.~Di~Francesco and G.~Stivaletta.
	\newblock Convergence of the follow-the-leader scheme for scalar conservation
	laws with space dependent flux.
	\newblock {\em Discrete Cont. Dyn. Syst.}, 40:233--266, 2020.
	
	\bibitem{ElKhatibGoatinRosini}
	N.~El-Khatib, P.~Goatin, and M.~D. Rosini.
	\newblock On entropy weak solutions of {H}ughes' model for pedestrian motion.
	\newblock {\em Z. Angew. Math. Phys.}, 64(2):223--251, 2013.
	
	\bibitem{GaravelloPiccolibook}
	M.~Garavello and B.~Piccoli.
	\newblock {\em Traffic flow on networks}, volume~1 of {\em AIMS Series on
		Applied Mathematics}.
	\newblock American Institute of Mathematical Sciences (AIMS), Springfield, MO,
	2006.
	\newblock Conservation laws models.
	
	\bibitem{MR4180806}
	L.~Gibelli, editor.
	\newblock {\em Crowd dynamics. {V}ol. 2--theory, models, and applications}.
	\newblock Modeling and Simulation in Science, Engineering and Technology.
	Birkh\"{a}user/Springer, Cham, 2020.
	
	\bibitem{MR3932134}
	L.~Gibelli and N.~Bellomo, editors.
	\newblock {\em Crowd dynamics. {V}ol. 1. {T}heory, models, and safety
		problems}.
	\newblock Modeling and Simulation in Science, Engineering and Technology.
	Birkh\"{a}user/Springer, Cham, 2018.
	
	\bibitem{bellomocrowd}
	L.~Gibelli and N.~Bellomo, editors.
	\newblock {\em Crowd Dynamics. Vol. 3}.
	\newblock Modeling and Simulation in Science, Engineering and Technology.
	Birkh\"{a}user/Springer, Cham, 2021.
	\newblock Modeling and Social Applications in the Time of COVID-19.
	
	\bibitem{MR3055243}
	P.~Goatin and M.~Mimault.
	\newblock The wave-front tracking algorithm for {H}ughes' model of pedestrian
	motion.
	\newblock {\em SIAM J. Sci. Comput.}, 35(3):B606--B622, 2013.
	
	\bibitem{MR4096596}
	M.~Gokieli and A.~Szczepa\'{n}czyk.
	\newblock A numerical scheme for evacuation dynamics.
	\newblock In {\em Parallel {P}rocessing and {A}pplied {M}athematics. {P}art
		{II}}, volume 12044 of {\em Lecture Notes in Comput. Sci.}, pages 277--286.
	Springer, Cham, [2020] \copyright 2020.
	
	\bibitem{HerzogPietschmannWinkler}
	R.~Herzog, J.-F. Pietschmann, and M.~Winkler.
	\newblock Optimal control of hughes' model for pedestrian flow via local
	attraction, 2020.
	
	\bibitem{HoldenRisebro02}
	H.~Holden and N.~H. Risebro.
	\newblock The continuum limit of {F}ollow-the-{L}eader models---a short proof.
	\newblock {\em Discrete Contin. Dyn. Syst.}, 38(2):715--722, 2018.
	
	\bibitem{HoldenRisebro01}
	H.~Holden and N.~H. Risebro.
	\newblock Follow-the-leader models can be viewed as a numerical approximation
	to the {L}ighthill-{W}hitham-{R}ichards model for traffic flow.
	\newblock {\em Netw. Heterog. Media}, 13(3):409--421, 2018.
	
	\bibitem{HUANG2009127}
	L.~Huang, S.~Wong, M.~Zhang, C.-W. Shu, and W.~H. Lam.
	\newblock Revisiting {H}ughes’ dynamic continuum model for pedestrian flow
	and the development of an efficient solution algorithm.
	\newblock {\em Transportation Research Part B: Methodological}, 43(1):127--141,
	2009.
	
	\bibitem{Hughes02}
	R.~L. Hughes.
	\newblock A continuum theory for the flow of pedestrians.
	\newblock {\em Transportation Research Part B: Methodological}, 36(6):507--535,
	2002.
	
	\bibitem{Kruzhkov}
	S.~N. Kruzhkov.
	\newblock First order quasilinear equations with several independent variables.
	\newblock {\em Mat. Sb. (N.S.)}, 81 (123):228--255, 1970.
	
	\bibitem{MR1927887}
	P.~G. LeFloch.
	\newblock {\em Hyperbolic systems of conservation laws}.
	\newblock Lectures in Mathematics ETH Z\"{u}rich. Birkh\"{a}user Verlag, Basel,
	2002.
	
	\bibitem{LWR1}
	M.~Lighthill and G.~Whitham.
	\newblock On kinematic waves. {II.} {A} theory of traffic flow on long crowded
	roads.
	\newblock In {\em Royal Society of London. Series A, Mathematical and Physical
		Sciences}, volume 229, pages 317--345, 1955.
	
	\bibitem{MR3451862}
	M.~Mimault.
	\newblock Scalar conservation law with discontinuity arising in pedestrian
	modeling.
	\newblock In {\em Congr\`es {SMAI} 2013}, volume~45 of {\em ESAIM Proc.
		Surveys}, pages 493--501. EDP Sci., 2014.
	
	\bibitem{Panov_traces2}
	E.~Y. Panov.
	\newblock Existence of strong traces for quasi-solutions of multidimensional
	conservation laws.
	\newblock {\em J. Hyperbolic Differ. Equ.}, 4(4):729--770, 2007.
	
	\bibitem{MR2592291}
	E.~Y. Panov.
	\newblock Existence and strong pre-compactness properties for entropy solutions
	of a first-order quasilinear equation with discontinuous flux.
	\newblock {\em Arch. Rat. Mech. Anal.}, 195:643--673, 2010.
	
	\bibitem{PietschmannStoetznerWinkler}
	J.-F. Pietschmann, A.~Stötzner, and M.~Winkler.
	\newblock Numerical investigation of agent controlled pedestrian dynamics using
	a structure preserving finite volume scheme, 2023.
	
	\bibitem{LWR2}
	P.~I. Richards.
	\newblock Shock waves on the highway.
	\newblock {\em Operations Research}, 4(1):42--51, 1956.
	
	\bibitem{Rosinibook}
	M.~D. Rosini.
	\newblock {\em Macroscopic models for vehicular flows and crowd dynamics:
		theory and applications}.
	\newblock Understanding Complex Systems. Springer, Heidelberg, 2013.
	
	\bibitem{MR3277564}
	M.~Twarogowska, P.~Goatin, and R.~Duvigneau.
	\newblock Macroscopic modeling and simulations of room evacuation.
	\newblock {\em Appl. Math. Model.}, 38(24):5781--5795, 2014.
	
	\bibitem{MR1869441}
	A.~Vasseur.
	\newblock Strong traces for solutions of multidimensional scalar conservation
	laws.
	\newblock {\em Arch. Ration. Mech. Anal.}, 160(3):181--193, 2001.
	
\end{thebibliography}
\end{document}